\documentclass[11pt,a4paper,final, dvipsnames]{amsart}

\input xy
\xyoption{all}

\DeclareMathAlphabet{\mathpzc}{OT1}{pzc}{m}{it}

\usepackage{amsfonts,amsmath,amssymb,indentfirst,mathrsfs,amscd}
\usepackage[mathscr]{eucal}
\usepackage{graphicx}
\usepackage{nicefrac}
\usepackage{amsthm}
\usepackage{tikz-cd}

\usepackage{color}
\usepackage[utf8]{inputenc}
\usepackage[spanish, english]{babel}
\usepackage{booktabs}
\usepackage{multirow}
\usepackage{verbatim}
\usepackage{mathtools}
\usepackage[symbol]{footmisc}

\renewcommand{\thefootnote}{\fnsymbol{footnote}}

\usepackage{calrsfs}
\DeclareMathAlphabet{ \mathcal}{OMS}{zplm}{m}{n}

\usepackage[all]{xy}
\usepackage{graphicx}
\usepackage{stmaryrd}
\usepackage{enumitem}

\usepackage{empheq}
\usepackage{color}
\usepackage{xcolor}
\usepackage{dsfont}

\usepackage[%pdftex,
            colorlinks=true,
            urlcolor=Blue,
            breaklinks=true,
            pdftitle={},
            pdfauthor={}]{hyperref}
\usepackage{color}
\definecolor{tocolor}{rgb}{.1,.1,.5}
\definecolor{urlcolor}{rgb}{.2,.2,.6}
\definecolor{linkcolor}{rgb}{.0,.3,.6}
\definecolor{citecolor}{rgb}{.6,.2,.2}
\hypersetup{colorlinks=true, urlcolor=urlcolor, linkcolor=linkcolor, citecolor=citecolor}

\author[]{\'Angel Gonz\'alez-Prieto}
\address{Facultad de Ciencias Matem\'aticas, Universidad Complutense  de  Madrid, Plaza Ciencias  3, 28040 Madrid Spain.}
\address{Instituto de Ciencias Matem\'aticas (CSIC-UAM-UC3M-UCM), C.\ Nicol\'as Cabrera 15, 28049 Madrid, Spain.}

\email{angelgonzalezprieto@ucm.es}

\title[]{Quantization of algebraic invariants through\\ Topological Quantum Field Theories}


\DeclareMathOperator{\coker}{coker\,}
\DeclareMathOperator{\Hom}{Hom\,}           %Hom%
\DeclareMathOperator{\tr}{Tr\,}             %Trace Tr%
\DeclareMathOperator{\colim}{colim}

\DeclareMathOperator{\Gr}{Gr}

\DeclareMathOperator{\Tr}{Tr\,}       %trace of a matrix%
\usepackage{amsmath}
\begin{document}

\newtheorem{thm}{Theorem}[section]
\newtheorem{prop}[thm]{Proposition}
\newtheorem{lem}[thm]{Lemma}
\newtheorem{cor}[thm]{Corollary}
\newtheorem{conjecture}{Conjecture}
\newtheorem*{theorem*}{Theorem}
\newtheorem*{question*}{Question}

\theoremstyle{definition}
\newtheorem{defn}[thm]{Definition}
\newtheorem{ex}[thm]{Example}
\newtheorem{as}{Assumption}

\theoremstyle{remark}
\newtheorem{rmk}[thm]{Remark}

\theoremstyle{remark}
\newtheorem*{prf}{Proof}

\newcommand{\iacute}{\'{\i}} %i con acento%
\newcommand{\norm}[1]{\lVert#1\rVert} %norma%

\newcommand{\lto}{\longrightarrow}
\newcommand{\hra}{\hookrightarrow}

\newcommand{\suchthat}{\;\;|\;\;}
\newcommand{\dbar}{\overline{\partial}}

\newcommand{\cA}{\mathcal{A}}
\newcommand{\cC}{\mathcal{C}}
\newcommand{\cD}{\mathcal{D}}
\newcommand{\cE}{\mathcal{E}}
\newcommand{\cF}{\mathcal{F}}
\newcommand{\cG}{\mathcal{G}} %Gauge group%
\newcommand{\cI}{\mathcal{I}} %Gauge group%
\newcommand{\cO}{\mathcal{O}} %Holomorphic functions sheaf%
\newcommand{\cM}{\mathcal{M}} %Moduli space% %moduli of parabolic bundles% %moduli of U(p,q) bundles%
\newcommand{\cN}{\mathcal{N}} %Space of minimal points of the Morse function%
\newcommand{\cP}{\mathcal{P}} %Moduli of K(D) pairs%
\newcommand{\cQ}{\mathcal{Q}} %Moduli of K(D) pairs%
\newcommand{\cS}{\mathcal{S}} %Moduli of solutions of Hitchin's equations, contructed by Konno%
\newcommand{\cU}{\mathcal{U}} %Moduli of stable U(p,q) parabolic Higgs bundles%
\newcommand{\cJ}{\mathcal{J}}
\newcommand{\cX}{\mathcal{X}}
\newcommand{\cT}{\mathcal{T}}
\newcommand{\cV}{\mathcal{V}}
\newcommand{\cW}{\mathcal{W}}
\newcommand{\cB}{\mathcal{B}}
\newcommand{\cR}{\mathcal{R}}
\newcommand{\cH}{\mathcal{H}}
\newcommand{\cZ}{\mathcal{Z}}
\newcommand{\D}{\bar{B}}
\newcommand{\Ss}{\mathcal{D}}

\newcommand{\Ker}{\textrm{Ker}\,}
\renewcommand{\coker}{\textrm{coker}\,}

\newcommand{\ext}{\mathrm{ext}} % an extension%
\newcommand{\x}{\times}

\newcommand{\mM}{\mathscr{M}} %Meromorphic function sheaf%

\newcommand{\CC}{\mathbb{C}} %Complex numbers%
\newcommand{\QQ}{\mathbb{Q}} %Rational numbers%
\newcommand{\FF}{\mathbb{F}} %Rational numbers%
\newcommand{\PP}{\mathbb{P}} %projective space%
\newcommand{\HH}{\mathbb{H}} %Hypercohomology, quaternions..%
\newcommand{\RR}{\mathbb{R}} %Real numbers%
\newcommand{\ZZ}{\mathbb{Z}} %Integer numbers%
\newcommand{\NN}{\mathbb{N}} %Natural numbers%
\newcommand{\DD}{\mathbb{D}} %Natural numbers%

\renewcommand{\lg}{\mathfrak{g}} %Lie algebra of G%
\newcommand{\lh}{\mathfrak{h}} %Lie algebra of H%
\newcommand{\lu}{\mathfrak{u}} %Lie algebra of U%
\newcommand{\la}{\mathfrak{a}} %Lie algebra of A%
\newcommand{\lb}{\mathfrak{b}} %Lie algebra of B%
\newcommand{\lm}{\mathfrak{m}} %Lie algebra of M%
\newcommand{\lgl}{\mathfrak{gl}} %Lie algebra of GL%
\newcommand{\lZ}{\mathfrak{Z}} %Almost-TQFT%
\newcommand\Unit[1]{\mathds{1}_{#1}}

\newcommand{\too}{\longrightarrow}
\newcommand{\imat}{\sqrt{-1}} %i%
\newcommand{\tinyclk}{{\scriptscriptstyle \Taschenuhr}} %Tiny Clock Symbol from ifsym package%
\newcommand\restr[2]{\left.#1\right|_{#2}}
\newcommand\rtorus[1]{{\mathbb{T}^{#1}}}
\newcommand\actPartial{\overline{\partial}}
\newcommand\handle[2]{\mathcal{A}^{#1}_{#2}}

\newcommand\pastingArea[2]{S^{#1-1} \times \bar{B}^{#2-#1}}
\newcommand\pastingAreaPlus[2]{S^{#1} \times \bar{B}^{#2-#1-1}}
\newcommand\coordvector[2]{\left.\frac{\partial}{\partial {#1}}\right|_{#2}}

% Categories
\newcommand\Sets{\textbf{Set}}
\newcommand\Cat{\textbf{Cat}}
\newcommand\Top{\textbf{Top}}
\newcommand\Topc{\textbf{Top}_c}
\newcommand\Toppc{\textbf{Topp}_c}
\newcommand\HTop{\textbf{HoTop}}
\newcommand\SkTop{\textbf{SkTop}}
\newcommand\TopHLC{\textbf{Top}_{hlc}}
\newcommand\TopS{\textbf{Top}_\star}
\newcommand\Diff{\textbf{Diff}}
\newcommand\Diffc{\textbf{Diff}_c}
\newcommand\GRing{\textbf{GRing}}

\newcommand\BaseBord[3]{\mathbf{Bd}_{{#1 #3}}^{#2}}
\newcommand\Bord[1]{\BaseBord{#1}{}{}}
\newcommand\NBord[2]{\mathbf{NBdp}_{{#1}}(#2)}
\newcommand\NBordnp[1]{\mathbf{NBdp}_{{#1}}}
\newcommand\BordC[1]{\BaseBord{#1}{clr}{}}
\newcommand\BordE[1]{\BaseBord{}{}{}(#1)}
\newcommand\EBord[2]{\BaseBord{#1}{#2}{}}
\newcommand\Bordo[1]{\BaseBord{#1}{or}{}}
\newcommand\Bordp[1]{\mathbf{Bdp}_{{#1}}}
\newcommand\Bordpar[2]{\mathbf{Bd}_{{#1}}(#2)}
\newcommand\Bordppar[2]{\mathbf{Bdp}_{{#1}}(#2)}

\newcommand\Tubo[1]{\mathbf{Tb}_{#1}^0}
\newcommand\Tub[1]{\mathbf{Tb}_{#1}}
\newcommand\ETub[2]{\mathbf{Tb}_{#1}^{#2}}
\newcommand\Tubp[1]{\mathbf{Tbp}_{#1}}
\newcommand\Tubpo[1]{\mathbf{Tbp}_{#1}^0}
\newcommand\Tubppar[2]{\mathbf{Tbp}_{#1}(#2)}

\newcommand\Embc{\textbf{Emb}_c}
\newcommand\EEmbc[1]{\textbf{Emb}_c^{#1}}
\newcommand\Embpc{\textbf{Embp}_c}
\newcommand\Embparc[1]{\textbf{Emb}_c(#1)}
\newcommand\Embpparc[1]{\textbf{Embp}_c(#1)}
\newcommand\OpEmb{\textbf{OEmb}}
\newcommand\OpEmbDiff{\textbf{OEmbDiff}}
\newcommand\OpEEmb[1]{\textbf{OEmb}^{#1}}
\newcommand\OpEEmbc[1]{\textbf{OEmb}^{#1}_c}
\newcommand\OpEmbc{\textbf{OEmb}_c}
\newcommand\OpEmbpc{\textbf{OEmbp}_c}
\newcommand\OpEmbp{\textbf{OEmbp}}
\newcommand\SkOpEmb{\textbf{SkOEmb}}

\newcommand\CDefc{\textbf{CDef}_c}
\newcommand\CEDefc[1]{\textbf{CDef}_c^{#1}}
\newcommand\CDefpc{\textbf{CDefp}_c}
\newcommand\CDefparc[1]{\textbf{CDef}_c(#1)}
\newcommand\CDefpparc[1]{\textbf{CDefp}_c(#1)}

\newcommand\CEmbc{\textbf{CEmb}_c}
\newcommand\Conc{\textbf{Con}_c}
\newcommand\Conpc{\textbf{Conp}_c}
\newcommand\CEmb{\textbf{CEmb}}
\newcommand\CEEmbc[1]{\textbf{CEmb}_c^{#1}}
\newcommand\CEmbpc{\textbf{CEmbp}_c}
\newcommand\CEmbparc[1]{\textbf{CEmb}_c(#1)}
\newcommand\CEmbpparc[1]{\textbf{CEmbp}_c(#1)}

\newcommand\cSp{\cS_p}
\newcommand\cSpar[1]{\cS_{#1}}

\newcommand\Diffpc{\textbf{Diff}_c}

\newcommand\PVar[1]{\mathbf{PVar}_{#1}}
\newcommand\PVarC{\mathbf{PVar}_{\CC}}
\newcommand\Sr[1]{\mathrm{S}{#1}}

\newcommand\Bordpo[1]{\mathbf{Bdp}_{{#1}}^{or}}
\newcommand\ClBordp[1]{\mathbf{l}\NBord{#1}{}{}}
\newcommand\CTub[3]{\mathbf{Tb}_{{#1 #3}}^{#2}}
\newcommand\CTubp[1]{\CTub{#1}{}{}}
\newcommand\CTubpp[1]{\mathbf{Tbp}_{#1}{}{}}
\newcommand\CTubppo[1]{\mathbf{Tbp}_{#1}^0}
\newcommand\CClose[3]{\mathbf{Cl}_{{#1 #3}}^{#2}}
\newcommand\CClosep[1]{\CClose{#1}{}{}}
\newcommand\Obj[1]{\mathrm{Obj}(#1)}
\newcommand\Mor[1]{\mathrm{Mor}(#1)}
\newcommand\Vect[1]{{#1}\textrm{-}\mathbf{Vect}}
\newcommand\Vecto[1]{{#1}\textrm{-}\mathbf{Vect}_0}
\newcommand\Mod[1]{{#1}\textrm{-}\mathbf{Mod}}
\newcommand\Modt[1]{{#1}\textrm{-}\mathbf{Mod}_t}
\newcommand\Rng{\mathbf{Ring}}
\newcommand\Grp{\mathbf{Grp}}
\newcommand\Grpd{\mathbf{Grpd}}
\newcommand\Grpdo{\mathbf{Grpd}_0}
\newcommand\HS[1]{\mathbf{HS}^{#1}}
\newcommand\MHS[1]{\mathbf{MHS}}
\newcommand\PHS[2]{\mathbf{HS}^{#1}_{#2}}
\newcommand\MHSq{\mathbf{HS}}
\newcommand\Sch{\mathbf{Sch}}
\newcommand\Sh[1]{\mathbf{Sh}\left(#1\right)}
\newcommand\QSh[1]{\mathbf{QSh}\left(#1\right)}
\newcommand\Var[1]{\mathbf{Var}_{#1}}
\newcommand\Varrel[1]{\mathbf{Var}/{#1}}
\newcommand\RVar[2]{\mathbf{Var}_{#1}/{#2}}
\newcommand\CVar{\Var{\CC}}
\newcommand\PHM[2]{\cM_{#1}^p(#2)}
\newcommand\MHM[1]{\cM_{#1}}
\newcommand\HM[2]{\textrm{HM}^{#1}(#2)}
\newcommand\HMW[1]{\textrm{HMW}(#1)}
\newcommand\VMHS[1]{VMHS({#1})}
\newcommand\geoVMHS[1]{VMHS_g({#1})}
\newcommand\goodVMHS[1]{\mathrm{VMHS}_0({#1})}
\newcommand\Par[1]{\mathrm{Par}({#1})}
\newcommand\K[1]{\mathrm{K}#1}
\newcommand\Ko[1]{\mathrm{K}{#1}_0}
\newcommand\KM[1]{\mathrm{K}\MHM{#1}}
\newcommand\KMo[1]{\mathrm{K}{\MHM{#1}}_0}
\newcommand\Ab{\mathbf{Ab}}
\newcommand\CPP{\cP\cP}
\newcommand\Bim[1]{{#1}\textrm{-}\mathbf{Bim}}
\newcommand\Span[1]{\mathbf{Span}({#1})}
\newcommand\ESpan[2]{\mathbf{Span}_{#2}({#1})}
\newcommand\Spano[1]{\mathbf{CoSpan}({#1})}
\newcommand\ESpano[2]{\mathbf{CoSpan}_{#2}({#1})}
\newcommand\GL[1]{\mathrm{GL}_{#1}}
\newcommand\SL[1]{\mathrm{SL}_{#1}}
\newcommand\PGL[1]{\mathrm{PGL}_{#1}}
\newcommand\AGL[1]{\mathrm{AGL}_{#1}}
\newcommand\Rep[1]{\mathfrak{X}_{#1}}
\newcommand\Repred[1]{\mathfrak{X}_{#1}^{r}}
\newcommand\Repirred[1]{\mathfrak{X}_{#1}^{ir}}
\newcommand\Charred[1]{\mathcal{R}_{#1}^{r}}
\newcommand\Charirred[1]{\mathcal{R}_{#1}^{ir}}
\newcommand\Repdiag[1]{\mathfrak{X}_{#1}^{ps}}
\newcommand\Repkappa[2]{\mathfrak{X}_{#1}^{#2}}
\newcommand\Reput[1]{\mathfrak{X}_{#1}^{ut}}

\newcommand\Qtm[1]{\mathcal{Q}_{#1}}
\newcommand\sQtm[1]{\mathcal{Q}_{#1}^0}
\newcommand\Fld[1]{\mathcal{F}_{#1}}
\newcommand\Id{\mathrm{Id}}

% Hodge theory
\newcommand\DelHod[1]{e\left(#1\right)}
\newcommand\RDelHod{e}
\newcommand\eVect{\mathcal{E}}
\newcommand\e[1]{\eVect\left(#1\right)}
\newcommand\intMor[2]{\int_{#1}\,#2}

% Representation Varieties
\newcommand\Dom[1]{\mathcal{D}_{#1}}

\newcommand\Xf[1]{{X}_{#1}}					% Free non-parabolic
\newcommand\Xs[1]{\mathfrak{X}_{#1}}							% Surface group non-parabolic

\newcommand\Xft[2]{\overline{{X}}_{#1, #2}} % Free tr=2
\newcommand\Xst[2]{\overline{\mathfrak{X}}_{#1, #2}}			% Surface group tr=2

\newcommand\Xfp[2]{X_{#1, #2}}			% Free J+
\newcommand\Xsp[2]{\mathfrak{X}_{#1, #2}}						% Surface group J+

\newcommand\Xfd[2]{X_{#1; #2}}			% Free diagonal
\newcommand\Xsd[2]{\mathfrak{X}_{#1; #2}}						% Surface group diagonal

\newcommand\Xfm[3]{\mathcal{X}_{#1, #2; #3}}		% Free mixed
\newcommand\Xsm[3]{X_{#1, #2; #3}}					% Surface group mixed

\newcommand\XD[1]{#1^{D}}
\newcommand\XDh[1]{#1^{\delta}}
\newcommand\XU[1]{#1^{UT}}
\newcommand\XP[1]{#1^{U}}
\newcommand\XPh[1]{#1^{\upsilon}}
\newcommand\XI[1]{#1^{\iota}}
\newcommand\XTilde[1]{#1^{\varrho}}
\newcommand\Xred[1]{#1^{r}}
\newcommand\Xirred[1]{#1^{ir}}

\newcommand\Char[1]{\cR_{#1}}
\newcommand\Chars[1]{\cR_{#1}}
\newcommand\CharW[1]{\mathscr{R}_{#1}}

% Derived Category
\newcommand\Ch[1]{\textrm{Ch}\,{#1}}
\newcommand\Chp[1]{\textrm{Ch}^+{#1}}
\newcommand\Chm[1]{\textrm{Ch}^-{#1}}
\newcommand\Chb[1]{\textrm{Ch}^b{#1}}
\newcommand\Der[1]{\textrm{D}{#1}}
\newcommand\Dp[1]{\textrm{D}^+{#1}}
\newcommand\Dm[1]{\textrm{D}^-{#1}}
\newcommand\Db[1]{\textrm{D}^b{#1}}

% TQFT
\newcommand\Gs{\cG}
\newcommand\Gq{\cG_q}
\newcommand\Gg{\cG_c}
\newcommand\Zs[1]{Z_{#1}}
\newcommand\Zg[1]{Z^{gm}_{#1}}
\newcommand\cZg[1]{\cZ^{gm}_{#1}}

% Miscelany
\newcommand\RM[2]{R\left(\left.#1\right|#2\right)}
\newcommand\RMc[3]{R_{#1}\left(\left.#2\right|#3\right)}
\newcommand\set[1]{\left\{#1\right\}}
\newcommand{\Stab}{\textrm{Stab}\,} %\Stabilizer%
\newcommand\EuChS{E}             %Trace Tr%
\newcommand\EuCh[1]{E\left(#1\right)}             %Trace Tr%

\newcommand{\Ann}{\textrm{Ann}\,}
\newcommand{\Rad}{\textrm{Rad}\,}  
\newcommand\supp[1]{\mathrm{supp}{(#1)}}
\newcommand\coh[1]{\left[H_c^\bullet\hspace{-0.05cm}\left(#1\right)\right]}
\newcommand\Kclass[1]{\left[#1\right]}
\newcommand\Bt{B_t}
\newcommand\Be{B_e}
\newcommand\re{\textrm{Re}\,}
\newcommand\imag{\textrm{Im}\,}
\newcommand\Kahc{\textbf{K\"ah}_c}

\newcommand\Cone[1]{\textrm{Cone}\left({#1}\right)}
\newcommand\conic{\textrm{C}}
\newcommand\op[1]{{#1}^{op}}

% Mixed Hodge modules
\newcommand\Ccs[1]{C_{cs}(#1)}
\newcommand\can{\textrm{can}}
\newcommand\var{\textrm{var}}
\newcommand\Perv[1]{\textrm{Perv}(#1)}
\newcommand\DR[1]{\textrm{DR}{#1}}
\renewcommand\Gr[2]{\textrm{Gr}_{#1}^{#2}\,}
\newcommand\ChV{\textrm{Ch}\,}
\newcommand\VerD{^{\textrm{Ve}}\DD}
\newcommand\DHol[1]{\textrm{D}^b_{\textrm{hol}}(\cD_{#1})}
\newcommand\RegHol[1]{\Mod{\cD_{#1}}_{\textrm{rh}}}
\newcommand\Drh[1]{\textrm{D}^b_{\textrm{rh}}(\cD_{#1})}
\newcommand\Dcs[2]{\textrm{D}^b_{\textrm{cs}}({#1}; {#2})}
\newcommand{\rat}{\mathrm{rat}}
\newcommand{\dmod}{\mathrm{Dmod}}

\hyphenation{mul-ti-pli-ci-ty}

\hyphenation{mo-du-li}

\begin{abstract}
In this paper we investigate the problem of constructing Topological Quantum Field Theories (TQFTs) to quantize algebraic invariants. We exhibit necessary conditions for quantizability based on Euler characteristics. In the case of surfaces, also provide a partial answer in terms of sufficient conditions by means of almost-TQFTs and almost-Frobenius algebras for wide TQFTs. As an application, we show that the Poincar\'e polynomial of $G$-representation varieties is not a quantizable invariant by means of a monoidal TQFTs for any algebraic group $G$ of positive dimension.
\end{abstract}
\null
\vspace{-1.1cm}
\maketitle

\vspace{-0.8cm}

%%%%%%%%%%%%%%%%%%%%%%%%%%%%%%%%%%%%%%%%%%%%%%%%%%%%%%%%%%%%%%%%
%%%%%%%%%%%%%%%%%%% SECTION: INTRODUCTION %%%%%%%%%%%%%%%%%%%%%%
%%%%%%%%%%%%%%%%%%%%%%%%%%%%%%%%%%%%%%%%%%%%%%%%%%%%%%%%%%%%%%%%

\section{Introduction}
{\let\thefootnote\relax\footnotetext{\noindent \hspace{-0.79cm}  \emph{2020 Mathematics Subject Classification}. Primary:
  57R56, % TQFT
 Secondary:
 18M05, % Monoidal categories, symmetric monoidal categories
 57K16, % Finite-type and quantum invariants, topological quantum field theories (TQFT)
 14D21. % Applications of moduli in mathematical physics

\noindent \emph{Key words and phrases}: Topological Quantum Field Theory, TQFT, quantization, monoidal structure, representation variety.}}

Since their inception by Witten \cite{witten1989quantum} and Atiyah \cite{Atiyah:1988}, Topological Quantum Field Theories (TQFTs) have become very important algebro-geometric tools in mathematical physics, homotopy theory, and knot theory, among others. Mathematically, a (monoidal) TQFT is a monoidal symmetric functor
$$
	\cZ: \Bord{n +1} \to \Mod{R}
$$
out of the category $\Bord{n +1} $ of $(n+1)$-dimensional bordisms into the category $\Mod{R}$ of modules over a fixed commutative unitary ring $R$, seen as a symmetric monoidal category.

Exploiting the monoidality of $\cZ$, the classification of such functors has been accomplished in the literature in the case of $(1+1)$-dimensional TQFTs \cite{Kock:2004}, fully extended TQFT as functors of $\infty$-categories \cite{lurie2009classification} or for extended $2$-dimensional TQFTs \cite{schommer2009classification} in terms of generators and relators of the bordism category. Indeed, extending these TQFTs to the $\infty$-category framework has motivated some of the exciting recent developments in higher category theory \cite{lurie2009higher}, motivic homotopy theory \cite{dundas2007motivic}, and factorization homology \cite{ayala2015factorization}. 

Apart from the theoretical developments, TQFTs have also been used to compute algebraic invariants of smooth manifolds. Indeed, suppose that we are interested in some invariant $\chi(W) \in R$ for $(n+1)$-dimensional orientable manifolds $W$, which may be quite hard to compute, e.g.\, because it captures the geometry of some moduli space attached to $W$. We will say that the invariant $\chi$ has been strongly quantized, if we manage to find a TQFT
$
	\cZ: \Bord{n +1} \to \Mod{R}
$
such that, for any closed $(n+1)$-dimensional manifold $W$, seen as a bordism $W: \emptyset \to \emptyset$, the map $\cZ(W): R \to R$ is given by multiplication by $\chi(W)$.

\begin{figure}[h!]
\includegraphics[width=0.6\textwidth]{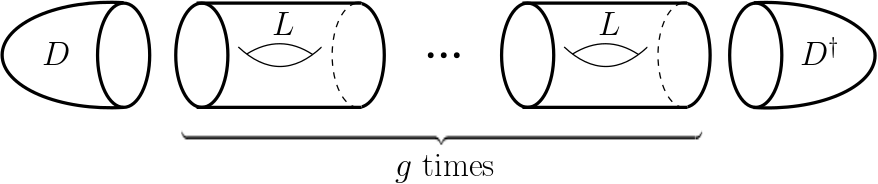}\vspace{-0.2cm}
\caption{Decomposition of a closed orientable surface into simpler pieces}
\label{fig:decomp-sigmag-intro}
\centering
\end{figure}

In this case, we can exploit the categorical structure of $\Bord{n+1}$ to decompose $W$ into simpler pieces and to re-ensemble the invariant by gluing the morphisms. For instance, in the $(1+1)$-dimensional case (surfaces) we have that the genus $g$ closed orientable surface $\Sigma_g$ can be decomposed as shown in Figure \ref{fig:decomp-sigmag-intro}.

Hence, in this context, we can compute the invariant $\chi(\Sigma_g)$ once at a time for all the surfaces as
\begin{equation}\label{eq:decomposition}
	\chi(\Sigma_g) = \cZ(D^\dag) \circ \left(\cZ(L)\right)^g \circ \cZ(D)(1).
\end{equation}

However, quantizing $\chi$ may be an impossible task for some important invariants. The reason is that the monoidality constraint of the TQFT actually imposes very strong restrictions, such as duality, on the modules associated through $\cZ$. The main aim of this paper is to investigate this quantizability problem: 

\begin{question*}
Given an invariant $\chi$, does there exist a TQFT $\cZ: \Bord{n +1} \to \Mod{R}$ such that $\chi(W) = \cZ(W)(1)$ for any closed orientable $(n+1)$-dimensional manifold?
\end{question*}

The first result in this direction that we show in this paper is the fact that if the invariant associated to the torus $\chi\left((S^1)^{n+1}\right)$ is not an integer of the ring $R$, then $\chi$ is not strongly quantizable. This is a direct consequence of a well-known argument regarding the trace of the map associated by $\cZ$ (see Proposition \ref{prop:necessity-quantization}).
However, despite its simplicity, from this result we directly get outstanding consequences. Suppose that we want to understand a certain family of moduli spaces $\cM_W$ attached to closed manifolds $W$ of dimension $n+1$. If $\cM_{(S^1)^{n+1}}$ is not an acyclic space, then the Poincar\'e polynomial of $\cM_W$ (or any other graded homology-based invariant of $\cM$) is not strongly quantizable.

For these reasons, the business of finding genuine TQFTs to compute homology invariants of classical moduli spaces is intrinsically too restrictive. To overcome this problem, in \cite{GPLM-2017}, it was proposed that a weak version, called an almost-TQFTs, may be used as substitute of genuine monoidal TQFTs. To be precise, let $\Tub{n+1}$ be the wide subcategory of $\Bord{n+1}$ with the same objects but whose morphisms are disjoint unions of `tubes', that is, bordisms with connected in and out boundaries. Then, an almost-TQFT is a monoidal symmetric functor
$$
	Z: \Tub{n+1} \to \Mod{R}.
$$
Of course, working in $\Tub{n+1}$ we are not allowed to decompose our manifold into any simple pieces, only into tubes. However, this restriction does not jeopardize most of the important gluing arguments used to compute invariants, such as Equation (\ref{eq:decomposition}), so almost-TQFTs can still be used to efficiently compute invariants. In this setting, we will say that an invariant $\chi$ is almost-quantizable if there exists an almost-TQFT that computes it. 

Thanks to the simple structure of $\Tub{n+1}$, it turns out that the almost-quantizability problem is much simpler than the strong quantizability. For instance, for surfaces,  $\Tub{1+1}$ is made of disjoint unions of the bordisms $D, D^\dag$ and $L$ of Figure \ref{fig:decomp-sigmag-intro}. In this way, in the case that the ground ring is a field $R = k$ and $Z$ takes values on finite dimensional $k$-vector spaces, an almost-TQFT for surfaces can be understood as a linear discrete dynamical system. Using that solutions to linear dynamical systems can be written explicitly, we get the following result.

\begin{theorem*}[Theorem \ref{thm:almost-quantizable} and Remark \ref{rmk:almost-quantizable-diagonal}]
A non-constant invariant $\chi$ of closed orientable surfaces with values in an algebraically closed field $k$ is almost-quantizable if and only if it has the form
\begin{equation*}
	\chi(\Sigma_g) = \sum_{i=1}^s \sum_{j=0}^{m_i-1} a_{i,j}\begin{pmatrix}g \\ j\end{pmatrix}\lambda_i^{g-j}.
\end{equation*}
for certain coefficients $a_{i,j}, \lambda_i \in k$ independent of $g$.
Furthermore, if the $N \times N$ matrix $\left(\chi(\Sigma_{g + h})\right)_{g,h}$ is invertible, where $N$ is the number of terms of the previous sum, then $\chi(\Sigma_g)$ for $g \geq 0$ must take the simpler form
\begin{equation*}
	\chi(\Sigma_g) = a_1\lambda_1^g + a_2\lambda_2^g + \ldots + a_N\lambda_N^g.
\end{equation*}
\end{theorem*}

This theorem is not only an existence result. As we show in Remark \ref{rmk:computation-tqft}, there is an algorithmic way of recovering the associated almost-TQFT by means of the knowledge of the invariants $\chi(\Sigma_g)$ for $g \leq 2N-1$. This observation turns the problem of computing such invariants, which is usually a very hard task that requires the use of advanced motivic or arithmetic techniques, into a problem of computing a finite set of examples that can be solved through a brute force approach.

With this method, the aforementioned result provides an effective criterion to decide whether an algebraic invariant is almost-quantizable. However, we still need a procedure to decide whether such almost-TQFT can be promoted to a genuine TQFT. To address this question, in this paper we provide the following partial result (for the definition of wide almost-TQFT, see Definition \ref{sec:monoidal-quant-almost-Frob}).

%in this paper we introduce almost-Frobenius algebras as the almost-TQFT counterpart of Frobenius algebras (which are known to be equivalent to $1+1$ monoidal TQFT \cite{Kock:2004}).

%Roughly speaking, an almost-Frobenius algebra is a tuple $(V, \tilde{T}, \tilde{\epsilon}, \tilde{\eta})$ with $V$ a finite dimensional $k$-vector space and linear maps $\tilde{T}: V \to V$, $\tilde{\epsilon}: k \to V$ and $\tilde{\eta}: V \to k$. These linear maps correspond to the images of the bordisms $L$, $D$ and $D^\dag$ of Figure \ref{fig:decomp-sigmag-intro} under the almost-TQFT. In this manner, the problem of promoting the almost-Frobenius algebra to a genuine Frobenius algebra is recasted as the existence of algebra $m: V \otimes_k V \to V$ and co-algebra $w: V \to V \otimes_k V$ structures on $V$ such that $\tilde{T} = m \circ w$. In general, there may not exist such a extension, or there can be infinitely many, but in the case that $V$ is spanned by the iterated vectors $v_g = \tilde{T}^g(\tilde{\epsilon}(1))$ for $g \geq 0$ (a situation that we shall call a wide almost-Frobenius algebra) it turns out that the problem is completely geometric and, at most, only one such extension may exist. Using these facts, we get the following criterion.

\begin{theorem*}[Theorem \ref{thm:quantization-almostTQFT}]
A $(1+1)$-dimensional wide almost-TQFT $Z$ with values in finite dimensional vector spaces can be extended to a TQFT if and only if the following conditions hold:
\begin{enumerate}
	\item\label{cond:1} The bilinear form $B(v_g, v_{g'}) = Z(D^\dag)(v_{g+g'})$ is non-degenerate.
	\item\label{cond:2} Let $b_1, \ldots, b_N$ be any orthogonal basis of $Z(S^1)$ with respect to $B$. Then, for all $1 \leq j \leq N$, we have
	$$
		Z(L)(b_j) = \sum_{i=1}^N \frac{b_jb_i^2}{B(b_i, b_i)},
	$$
	where $b_ib_j = \sum_{g,g'} a_i^g a^{g'}_j\, v_{g+g'}$ for $b_i = \sum_g a^g_i\, v_{g}$ and $b_j = \sum_{g'} a^{g'}_j\, v_{g'}$.
\end{enumerate}
\end{theorem*}

As an application of these results, we shall focus on representation varieties. This is the space parametrizing, for a closed manifold $W$, the collection of representations of the fundamental group of $W$ into a fixed complex algebraic group $G$
$$
	\Rep{G}(W) = \Hom(\pi_1(W), G).
$$
By choosing a presentation of $\pi_1(W)$, the set $\Rep{G}(W)$ can be naturally endowed with the structure of an algebraic variety, called the $G$-representation variety. In the case $G = \GL{n}(\CC), \SL{n}(\CC)$ and $W$ a closed orientable surface, these spaces play a crucial role in the non-abelian Hodge correspondence since their GIT quotients under the adjoint action (the so-called character varieties) turn out to be homeomorphic to moduli spaces of Higgs bundles \cite{Corlette:1988} and of flat connections \cite{SimpsonI,SimpsonII}.

A direct dimensional count shows that $\Rep{G}(S^1 \times S^1)$ is not an acyclic space for $\dim G > 0$ (Theorem \ref{thm:non-quantizable-rep-var}), so a direct application of the simple criterion of Corollary \ref{cor:non-quantizable-integers} shows that graded homology invariants of the representation variety are not quantizable. A posteriori, this behavior is clearly expectable since TQFTs do not take into account basepoints, which are crucial for the functoriality of the fundamental group.

For this reason, any hope of quantizing homology invariants of representation varieties requieres to enlarge the category of bordism to equip its objects with a finite number of basepoints, leading to the category $\Bordp{n+1}$ of bordisms with baspoints. This extension also gives rise to the so-called (monoidal) TQFTs with basepoints, which are monoidal symmetric functors $\cZ: \Bordp{n+1} \to \Mod{R}$ (see Section \ref{sec:tqfts-with-basepoints}).

It turns out that the simple argument used in the basepoint-free case to show that invariants of tori must be integers no longer works in this setting. In some sense, TQFTs with basepoints are more flexible and give more room to quantize invariants, since we are adding some extra information to the bordisms. However, by studying the behaviour of TQFTs on trivial bordisms with two basepoints, one for each boundary component, analogous necessary conditions for strong quantizability with basepoints can be provided (see Proposition \ref{prop:trace-twisting}). From this restriction and a subtler analysis on the motivic structure of the representation variety, we get the following result.

\begin{theorem*}[Theorems \ref{thm:non-quantizable-rep-var}, \ref{thm:non-quantizable-basepoints-rep-var} and \ref{thm:lax-monoidal-rep-var}]
Let $G$ be a complex algebraic group of dimension $\dim G > 0$. Then, the Poincar\'e polynomial of $G$-representation varieties of closed orientable manifolds of dimension $\geq 2$:
\begin{enumerate}
	\item Is not strongly quantizable.
	\item Is not strongly quantizable with basepoints (and split).
	\item Is almost-quantizable with basepoints.
\end{enumerate}
\end{theorem*}

Of course, the same result can be applied to more refined invariants of the complex structure on the representation variety, such as its Hodge polynomial, its $E$-polynomial or even its virtual class (motive) in the Grothendieck ring of algebraic varieties. The no-go statements (1) and (2) are obtained in this paper for the first time, settling an important open question conjectured in \cite{MM}. Statement (3) was obtained in \cite{GP-2018}, but the lax monoidality of the built TQFT was though to be a limitation of the construction, and it was open whether such weak TQFT could be promoted to a monoidal TQFT. The results of this paper show that the construction of \cite{GP-2018} is actually sharp.

On the other hand, in the case that $G$ is a finite group (in particular if $\dim G = 0$), the number of points of the $G$-representation variety is obviously an integer, so the necessary conditions of Corollary \ref{cor:non-quantizable-integers} are fulfilled. In fact, in this situation, it turns out that the almost-TQFT constructed in \cite{GP-2018} can be promoted to a genuine TQFT, giving rise to the following result.

\begin{theorem*}[\cite{GP-2018} and Proposition \ref{prop:TQFT-dim0-monoidal}]
Let $G$ be a finite group. Then, the number of closed points of $G$-representation varieties is strongly quantizable with basepoints.
\end{theorem*}

Finally, as an application of the interpretation of a $(1+1)$-dimensional almost-TQFT as a dynamical system, in this paper we also provide a novel very simple method to compute the $E$-polynomial of representation varieties over arbitrary orientable surfaces from the knowledge of the $E$-polynomial in a finite number of cases (which can be obtained through brute force). The approach is somehow heuristic, since it relies on the fact that the almost-TQFT computing virtual Hodge classes of $G$-representation varieties restricts to a finite dimensional almost-TQFT, a question that currently remains open. This fact is known in the case $G = \SL{2}(\CC)$, so only the knowledge of the virtual Hodge classes for surfaces of genus $0 \leq g \leq 11$ is needed and from this we can re-prove the following result of \cite[Proposition 11]{MM} (see \cite[Remark 5.11]{GP-2018} for the case of virtual classes).

\begin{theorem*}
The $E$-polynomial of the $\SL{2}(\CC)$-representation variety is
\begin{align*}
	[\Rep{\SL{2}(\CC)}(\Sigma_g)] & = \,{\left(q^2 - 1\right)}^{2g - 1} q^{2g - 1} +
\frac{1}{2} \, {\left(q -
1\right)}^{2g - 1}q^{2g -
1}(q+1){\left({2^{2g} + q - 3}\right)} 
\\ & \qquad+ \frac{1}{2} \,
{\left(q + 1\right)}^{2g + r - 1} q^{2g - 1} (q-1){\left({2^{2g} +q -1}\right)} + q(q^2-1)^{2g-1},
\end{align*}
where $q = uv = e(\mathbb{A}^1_\CC)$ is the $E$-polynomial of the affine line.
\end{theorem*}

It is worth mentioning that, in this work, we have focused on representation varieties since they are the building blocks of other moduli spaces such as character varieties (through GIT quotient) or character stacks (through stacky quotients). Nevertheless, the same type of arguments shows that, for $G = \SL{2}(\CC)$, neither character varieties (from the computation of its $E$-polynomial in \cite{MM}) nor character stacks (from the computation of its virtual class in \cite{gonzalez2022virtual}) are acyclic, so their homology invariants are not strongly quantizable either. In the case of character stacks, they turn out to be almost-quantizable \cite{gonzalez2022virtual}, but almost nothing is known for character varieties (apart from the extended character field theory of Ben-Zvi, Gunningham and Nadler \cite{Ben-Zvi-Gunningham-Nadler}). From these results, we may expect that character stacks are not strongly quantizable and character varieties are not even almost-quantizable.

\subsection*{Structure of the manuscript}
Section \ref{sec:topological-rest} reviews the fundamental properties of TQFTs, discussing the quantizability problem (Section \ref{sec:quantizability-problem}) and studying almost-TQFTs and lax monoidal TQFTs (Section \ref{sec:almost-tqfts}). Necessary conditions for monoidality are given in that section, as well as the equivalence between lax monoidal and almost-TQFTs. Section \ref{sec:almost-frob} is devoted to sufficient conditions for monoidality in the case of $(1+1)$-dimensional TQFTs (surfaces). In particular, Section \ref{sec:monoidal-quant-almost-Frob} studies almost-Frobenius algebras and their interplay with usual Frobenius algebras, whereas Section \ref{sec:dynamical-system} discusses the interpretation of almost-TQFTs as linear discrete dynamical systems. The case of representation varieties is analyzed in Section \ref{sec:representation-varieties}, where TQFTs with basepoints are proposed in Section \ref{sec:tqfts-with-basepoints} and the corresponding quantizability problem is studied in Sections \ref{sec:split-tqfts-rep} and \ref{sec:lax-monoidal-rep-var}. Finally, Section \ref{sec:sl2-rep-var} exemplifies the new computational method for $E$-polynomials based on the interpretation of the almost-TQFT as a dynamical system for $\SL{2}(\CC)$-representation varieties.

\subsection*{Acknowledgments}

The author is greatly indebted to Thomas Wasserman for very enlightening conversations regarding the monoidality problem. He is also very grateful to Jesse Vogel for his very careful reading of the manuscript, his invaluable and insightful comments to improve the paper, and for pointing out several mistakes in the first version of this manuscript. The author also wants to thank M. Ballandras, G. Barajas, G. Gallego, M. Hablicsek, K. Knop and A. Saha for useful discussions.

The author also acknowledges the hospitality of Department of Mathematics at Universidad Aut\'onoma de Madrid where this work was partially completed. This work has been partially supported by the \textit{Madrid Government (Comunidad de Madrid -- Spain)} under the Multiannual Agreement with the Universidad Complutense de Madrid in the line Research Incentive for Young PhDs, in the context of the V PRICIT (Regional Programme of Research and Technological Innovation) through the project PR27/21-029, the Ministerio de Ciencia e Innovaci\'on Project PID2021-124440NB-I00 (Spain) and the BBVA Foundation COMPLEXFLUIDS project.

%%%%%%%%%%%%%%%%%%%%%%%%%%%%%%%%%%%%%%%%%%%%%%%%%%%%%%%%%%%%%%%%
%%%%%%%%%%%%%%%%%%%%%%% BIBLIOGRAPHY %%%%%%%%%%%%%%%%%%%%%%%%%%%
%%%%%%%%%%%%%%%%%%%%%%%%%%%%%%%%%%%%%%%%%%%%%%%%%%%%%%%%%%%%%%%%

\section{Topological restrictions to the monoidality problem}\label{sec:topological-rest}

A monoidal category is a tuple $(\cC, \otimes, I)$ comprised of a category $\cC$ equipped with a bifunctor $\otimes: \cC \times \cC \to \cC$ and an object $I$ of $\cC$, usually referred to as the unit of the monoidal structure, such that we have natural equivalences $- \otimes c \cong c \otimes - \cong \Id_c$ for any object $c$ of $\cC$. Some extra coherent conditions are typically required, as described in \cite[Section VII.1]{MacLane}.

For our purposes, the most important example of monoidal category is the following. Let $n \geq 0$ be a natural number. The category of $(n+1)$-dimensional bordisms $\Bord{n+1}$ is the category comprised of the following information:
\begin{itemize}
	\item Objects: The objects of $\Bord{n+1}$ are the closed (i.e.\ compact and boundaryless) orientable $n$-dimensional smooth manifolds, possibly empty.
	\item Morphisms: A morphism $M_1 \to M_2$ is a class of $(n+1)$-dimensional orientable bordisms, that is, compact smooth orientable manifolds $W$ such that $\partial W=M_1 \sqcup M_2$. Two such bordisms are declared as equivalent if there exists a boundary-preserving diffeomorphism between them.
	\item Composition: Given two bordisms $W: M_1 \to M_2$ and $W': M_2 \to M_3$, the composed morphism $W' \circ W: M_1 \to M_3$ is the result of gluing $W$ and $W'$ along their common boundary $M_2$, usually denoted by $W \cup_{M_2} W'$. 
	\item Monoidality: $\Bord{n+1}$ is endowed with its natural monoidal structure given by disjoint union (of both closed manifolds and bordisms). The unit of the monoidal structure is the empty manifold $\emptyset$.
\end{itemize}

\begin{rmk}
Notice that differentiable gluing is only well-defined up to difeomorphism. That is why we should take the morphisms as classes of cobordisms up to boundary-preserving diffeomorphism.
\end{rmk}

Given monoidal categories $(\cC, \otimes_\cC, I_\cC)$ and $(\cD, \otimes_\cD, I_\cD)$, a functor $\cF: \cC \to \cD$ is said to be \emph{lax monoidal} if there exists a family of morphisms in $\cD$
\begin{equation}\label{eq:lax-monoidal-maps}
	\Delta_{c,c'}: \cF(c \otimes_{\cC} c') \to \cF(c) \otimes_\cD \cF(c'), 
\end{equation}
natural in $c, c' \in \cC$ and commuting with the map interchanging the factors in the monoidal structures, as well as an isomorphism $\cF(I_\cC) \cong I_\cD$. The functor is said to be \emph{monoidal} if all the morphisms $\Delta_{c,c'}$ are isomorphisms.

\begin{defn}
Let $\cC$ a symmetric monoidal category. An \emph{$(n + 1)$-lax monoidal Topological Quantum Field Theory} (TQFT) on $\Bord{n + 1}$ with values in $\cC$ is a lax monoidal and symmetric functor $\cZ: \Bord{n + 1} \to \cC$. Furthermore, if $\cZ$ is a monoidal functor we will say that $\cZ$ is a \emph{monoidal TQFT}.
\end{defn}

\begin{rmk}
Typically, as in \cite{Atiyah:1988}, the target category is taken as $\cC = \Mod{R}$, the category of $R$-modules and module homomorphism between them, for a fixed commutative unitary ring $R$.
\end{rmk}

\begin{rmk}
In $\Bord{n+1}$ we took our objects and morphisms to be orientable manifolds (but without a fixed orientation). We make this choice to recover the classical result that shows that $(1+1)$-dimensional TQFTs are equivalent to Frobenius algebras \cite{Kock:2004} since otherwise we would need to add an extra non-orientable morphism. However, the orientation (and the smooth structure) plays no role in the following arguments and, with the appropriate modifications, most of the results of this paper can be adapted to work in the non-orientable setting. 
\end{rmk}

\subsection{The quantizability problem}\label{sec:quantizability-problem}

Let us consider an invariant $\chi$ of $(n+1)$-dimensional closed orientable manifolds with values in a set $\Lambda$ which is invariant under diffeomorphisms. This means that $\chi$ is a map that assigns, to any $(n+1)$-dimensional closed orientable manifold $W$, an element $\chi(W) \in \Lambda$. This element must be invariant under diffeomorphisms, that is, if $W$ and $W'$ are diffeomorphic then $\chi(W) = \chi(W')$.

Recall that a category $\cC$ is said to be \emph{locally small} if for any objects $c,c'$ of $\cC$, the morphisms $\Hom_\cC(c,c')$ is a set (instead of a proper class). If, in addition, the objects of $\cC$ form a set, then $\cC$ is said to be \emph{small}.

\begin{defn}\label{defn:quantizability}
Let $(\cC, \otimes, I)$ be a locally small symmetric monoidal category satisfying $\Lambda \subseteq \Hom_{\cC}(I, I)$. A $\Lambda$-valued diffeomorphism invariant $\chi$ is said to be \emph{lax-quantizable} in $\cC$ if there exists a lax monoidal TQFT
$$
	\cZ: \Bord{n+1} \to \cC
$$
such that for any $(n+1)$-dimensional closed orientable manifold $W$, seen as a bordism $W: \emptyset \to \emptyset$, we have that $Z(W) = \chi(W)$.

If $\cZ$ can be taken to be a monoidal TQFT, we shall say that $\chi$ is \emph{strongly quantizable}.
\end{defn}

\begin{rmk}
If $\cC = \Mod{R}$, then we have a natural identification $R \cong \Hom_{\Mod{R}}(R, R)$. In this case, a lax monoidal TQFT $\cZ: \Bord{n+1} \to \Mod{R}$ quantizes $\chi$ if $Z(W)(1) = \chi(W)$ for every $(n+1)$-dimensional closed orientable manifold $W$.
\end{rmk}

Recall that, given a symmetric monoidal category $(\cC, \otimes, I)$ and an object $c \in \cC$, a \emph{dual} of $c$ is an object $c^* \in \cC$ together with morphisms
$$
	\mu: c \otimes c^* \to I, \qquad \delta: I \to c^* \otimes c,
$$
such that $(\mu \otimes \Id_{c}) \circ (\Id_c \otimes \delta) = \Id_c$ and $(\Id_{c^*} \otimes \mu ) \circ ( \delta \otimes \Id_{c^*}) = \Id_{c^*}$. Dual objects are unique up to isomorphism. An object that admits a dual is said to be \textit{dualizable}. Given $f: c \to c$ with $c$ a dualizable object, the \emph{trace} of $f$ is the morphism $\Tr(f): I \to I$ given by the composition
$$
	I \stackrel{\delta}{\longrightarrow} c \otimes c^* \cong c^* \otimes c \stackrel{\Id_{c^*} \otimes f}{\longrightarrow} c^* \otimes c \stackrel{\mu}{\longrightarrow} I.
$$
The \emph{Euler characteristic} of a dualizable object $c \in \cC$ is $[c] = \Tr(\Id_{c}) = \mu \circ \delta$. For further information, see \cite{dold1983duality}.

It turns out that, for some manifolds, the invariant computed by a monoidal TQFT is just an Euler characteristic. The proof is well-known, but we include it here for completeness.

\begin{prop}\label{prop:necessity-quantization}
Let $\cZ: \Bord{n+1} \to \cC$ be a monoidal TQFT. Then, for any $n$-dimensional closed orientable manifold $M$ we have that the induced homomorphism
$$
	\cZ(M \times S^1): I \to I
$$ 
is the Euler characteristic of  $\cZ(M)$. 

\begin{proof}
Let us denote by $\delta = M \times [0,1]: \emptyset \to M \sqcup M$ and by $\mu = M \times [0,1]: M \sqcup M \to \emptyset$ the two `elbow' bordisms. As depicted in Figure \ref{fig:zorro}, the ``Zorro'' bordism shows that $M$ is self-dual with duality maps $\delta$ and $\mu$.
\begin{figure}[h!]
\includegraphics[width=9cm]{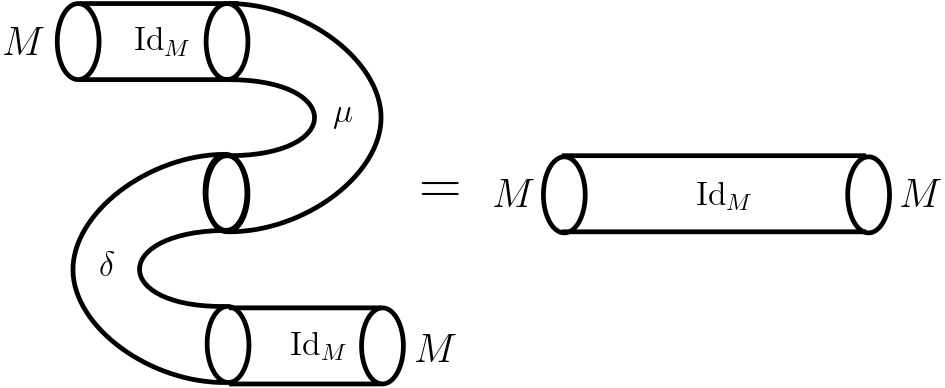}
\caption{The ``Zorro'' bordism}
\label{fig:zorro}
\centering
\end{figure}

Since a monoidal functor preserves duality, we have that $\cZ(M)$ is self-dual with duality maps $\cZ(\delta)$ and $\cZ(\mu)$. In particular, we have that
$$
	\cZ(M \times S^1) = \cZ(\mu \circ \delta) = \cZ(\mu) \circ \cZ(\delta) = \Tr(\Id_{\cZ(M)}).
$$
\end{proof}
\end{prop}

Using Proposition \ref{prop:necessity-quantization} with $\cC = \Mod{R}$, and noting that the Euler characteristic of a finitely generated projective $R$-module $N$ (the dualizable objects of $\Mod{R}$) is the rank of $N$, we find a useful criterion for non-quantizability of invariants. Recall that the subring of integers of a ring $R$ is the subring of $R$ generated by $1$.

\begin{cor}\label{cor:non-quantizable-integers}
Let $\chi$ be an $R$-valued invariant of $(n+1)$-dimensional closed orientable manifolds. If there exists an $n$-dimensional closed orientable manifold $M$ such that $\chi(M \times S^1)$ does not lie in the
subring of integers of $R$, then $\chi$ is not strongly quantizable. 
\end{cor}

\subsection{Almost-TQFTs}\label{sec:almost-tqfts} 

\begin{figure}[h!]
\includegraphics[width=0.75\textwidth]{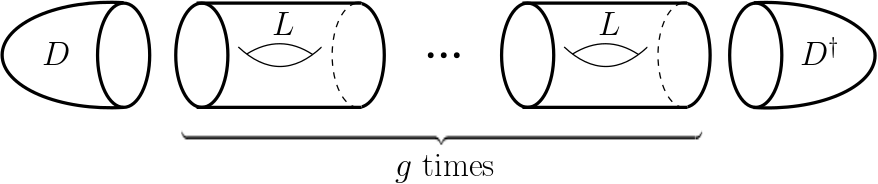}\vspace{-0.2cm}
\caption{Decomposition of a closed orientable surface into simpler pieces}
\label{fig:decomp-sigmag}
\centering
\end{figure}

Corollary \ref{cor:non-quantizable-integers} prevents many invariants to be strongly quantizable by means of a monoidal TQFT. However, with a view towards applications, full monoidality may be not needed. For instance, regarding $(1+1)$-dimensional TQFTs, any closed orientable surface can be decomposed into the composition of a disc $D: \emptyset \to S^1$, a bunch of holed torus $L: S^1 \to S^1$ and a disc in the other way around $D^\dag: S^1 \to \emptyset$, as shown in Figure \ref{fig:decomp-sigmag}. 

In this manner, to compute invariants we do not need to know the value of a TQFT on an arbitrary disjoint union of circles, but only on connected manifolds. This is a weaker class of functors that we shall call almost-TQFTs.

\begin{defn}
Fix an integer $n \geq 0$. The \emph{category $\Tubo{n + 1}$ of $(n+1)$-strict tubes} is the subcategory of $\Bord{n+1}$ whose objects and morphisms are connected manifolds (maybe empty). The \emph{category $\Tub{n+1}$ of $(n+1)$-tubes} is the monoidal closure of $\Tubo{n+1}$, that is, $\Tub{n+1}$ has the same objects as $\Bord{n+1}$ and a morphism $W$ of $\Bord{n+1}$ is in $\Tub{n+1}$ if $W$ can be decomposed as a disjoint union $W = \sqcup W_i$ of bordisms, where each $W_i$ has connected in and out boundaries.
\end{defn}

\begin{defn}
Let $(\cC, \otimes, I)$ be a symmetric monoidal category. An \emph{$(n+1)$-almost Topological Quantum Field Theory} with values in $\cC$ is a symmetric monoidal functor
$$
	Z: \Tub{n+1} \to \cC.$$
\end{defn}

In analogy with Definition \ref{defn:quantizability}, a $\Lambda$-valued diffeomorphism invariant $\chi$ is said to be \emph{almost-quantizable} if there exists an almost-TQFT
$ \cZ $
such that $\Lambda \subseteq \Hom_{\cC}(I,I)$ and $Z(W) = \chi(W)$ for any $(n+1)$-dimensional closed orientable manifold $W$.

Notice that any lax monoidal TQFT $\cZ: \Bord{n+1} \to \cC$ gives rise to a functor $\Tubo{n+1} \to \cC$ just by restriction. This functor is automatically promoted to an almost-TQFT $Z: \Tub{n+1} \to \cC$ by setting $Z(\sqcup_i M_i) = \otimes_i Z(M_i)$, both for objects and morphisms.  Hence, any lax monoidal TQFT induces an almost-TQFT. The following results gives a reciprocal to this observation.

Recall that a category is called \emph{(co)complete} if all the small (co)limits exist, i.e.\ limits over small categories. A functor is said to be \emph{(co)continuous} if it commutes with all the small (co)limits. An important concept for the upcoming argument is the one of a left Kan extension. Given functors $\iota: \cT \to \cB$ and $F: \cT \to \cC$, a \emph{left Kan extension} of $F$ along $\iota$ is a functor $\cF: \cB \to \cC$ with a natural transformation $\eta: F \Rightarrow \cF \circ \iota$. When $\cT$ is small and $\cC$ is cocomplete, left Kan extensions always exist. For more information, please refer to \cite[Chapter X]{MacLane}.

\begin{thm}\label{thm:kan-extension}
Suppose that the category $\cC$ is cocomplete and that $\otimes: \cC \times \cC \to \cC$ is bi-cocontinuous. Given an almost-TQFT $Z: \Tub{n+1} \to \cC$, there exists a lax monoidal TQFT $\cZ: \Bord{n+1} \to \cC$ extending $Z$.

\begin{proof}
The proof follows similar lines to \cite{fritz2018criterion}, but weaker hypotheses are needed due to the particular situation considered. Let us spell out the details. 

We take $\cZ$ as the left Kan extension of $Z$ along the inclusion $\iota: \Tub{n+1} \hookrightarrow \Bord{n+1}$
\[\xymatrix{
\Bord{n+1} \ar@{--{>}}[r]^{\cZ} & \cC \\
\Tub{n+1} \ar@{^{(}-{>}}[u]^{\iota} \ar[ru]_{Z} &
}
\]
This diagram is equipped with a natural transformation $\eta: Z \Rightarrow \cZ \circ \iota$.

Observe that this Kan functor is actually an extension for the objects of $\Tubo{n+1}$, i.e. $\eta_{M}$ is an isomorphism for $M \in \Tubo{n+1}$. This can be obtained through an adaptation of the argument of \cite[Corollary X.4]{MacLane}. To be precise, for $M \in \Tubo{n+1}$, the left Kan extension $\cZ(M)$ is given as the colimit of the functor $Z \circ Q_{M}: (\iota \downarrow \iota(M)) \to \cC$, where $Q_M: (\iota \downarrow \iota(M)) \to \Tubo{n+1}$ is $Q_M(M', f) = M'$ for an object $f: \iota(M') \to \iota(M)$ of the comma category $(\iota \downarrow \iota(M))$. Since the inclusion $\iota: \Tubo{n+1} \hookrightarrow \Bord{n+1}$ is full and faithful, then $\Id_M: M \to M$ is a final object of $(\iota \downarrow \iota(M))$. Hence, the colimit is computed by evaluating $Z \circ Q_{M}$ in this final object so $\cZ(M) = Z \circ Q_{M}(\Id_M, M) = Z(M)$.

Now, let us construct the lax monoidal structure on $\cZ$. Let us fix $M_1, M_2 \in \Tub{n+1}$. Consider the functor
$$
	(Z \circ Q_{M_1}) \otimes (Z \circ Q_{M_2}): (\iota \downarrow M_1) \times (\iota \downarrow M_2) \longrightarrow \cC \times \cC \stackrel{\otimes}{\longrightarrow} \cC.
$$
Since $\otimes$ is a bi-cocontinuous functor, we have that
$$
	\colim \left( (Z \circ Q_{M_1}) \otimes (Z \circ Q_{M_2})\right) = \colim (Z \circ Q_{M_1}) \otimes \colim (Z \circ Q_{M_2}) = \cZ(M_1) \otimes \cZ(M_2).
$$

By juxtaposing the bordisms, we get a natural inclusion $j: (\iota \downarrow M_1) \times (\iota \downarrow M_2) \to (\iota \downarrow M_1 \sqcup M_2)$. Hence, $\cZ(M_1 \sqcup M_2)$, which is the colimit in $(\iota \downarrow M_1 \sqcup M_2)$ of $Z \circ Q_{M_1 \sqcup M_2}$, is a cocone in $(\iota \downarrow M_1) \times (\iota \downarrow M_2)$ for $Z \circ Q_{M_1 \sqcup M_2} \circ j = (Z \circ Q_{M_1}) \otimes (Z \circ Q_{M_2})$. Therefore, we get a morphism
$$
	\Delta_{M_1, M_2}: \cZ(M_1) \otimes \cZ(M_2) \to \cZ(M_1 \sqcup M_2).
$$
A straightforward computation shows that the maps $\Delta_{M_1,M_2}$ satisfy all the properties needed for the lax monoidal structure maps. For further details, see \cite[Theorem 2.1]{fritz2018criterion}.
\end{proof}
\end{thm}

\begin{rmk}
In the case that $M_1$ and $M_2$ are connected, the lax monoidal map $\Delta_{M_1, M_2}: \cZ(M_1) \otimes_{R} \cZ(M_2) \to \cZ(M_1 \sqcup M_2)$ can be obtained from the natural transformation $\eta: Z \Rightarrow \cZ \circ \iota$. Indeed, since $\eta_{M_1}, \eta_{M_2}$ are isomorphisms, we can complete the diagram
\[\xymatrix{
Z(M_1) \otimes Z(M_2) = Z(M_1 \sqcup M_2) \ar[rr]^{\hspace{1.7cm} \eta_{M_1 \sqcup M_2}} \ar[d]_{\eta_{M_1} \otimes \eta_{M_2}} && \cZ(M_1 \sqcup M_2) \\
\cZ(M_1) \otimes \cZ(M_2) \ar@{--{>}}[rru]
}
\]
A direct computation shows that this map is actually a map of cocones. Hence, since the map from the colimit to a cocone is unique, this map $ \cZ(M_1) \otimes_{R} \cZ(M_2) \to \cZ(M_1 \sqcup M_2)$ coincides with the one of Theorem \ref{thm:kan-extension}.
\end{rmk}

\begin{ex}
The category $\cC = \Mod{R}$ of $R$-modules satisfies the hypotheses of Theorem \ref{thm:kan-extension}: It is cocomplete by \cite[Section V.1]{MacLane} and $\otimes$ is bi-cocontinuous since $- \otimes M$ has $\Hom(M, -)$ as right adjoint.
\end{ex}

\begin{cor}
An invariant is lax-quantizable if and only if it is almost-quantizable.
\end{cor}

\section{Almost-Frobenius algebras and sufficient conditions for monoidality}\label{sec:almost-frob}

A \emph{(commutative) Frobenius $k$-algebra} is a finite commutative $k$-algebra $A$ endowed with a non-degenerate symmetric bilinear form $B: A \times A \to k$ such that $B(xy,z) = B(x,yz)$ for all $x,y,z \in A$. It is well-known that the category of commutative Frobenius algebras (with intertwining maps between them as morphisms) is equivalent to the category of $(1+1)$-dimensional monoidal TQFTs with values in the category $\Vecto{k}$ of finite dimensional $k$-vector spaces \cite{Kock:2004}, where the ring structure on $A= \cZ(S^1)$ comes from the image of the pair of pants $S^1 \sqcup S^1 \to S^1$, the unit of the ring is the image of a disc, and the bilinear form $B$ corresponds to the `elbow' bordism $S^1 \sqcup S^1 \to \emptyset$.

In particular, several new maps can be computed from this information. The most relevant ones for our purposes are the linear maps
$$
	T = m \circ w: A \to A, \quad \epsilon: k \to A, \quad \eta: A \to k,
$$
where $m: A \otimes_k A \to A$ is the ring multiplication map, $w: A \to A \otimes_k A$ is the associated co-multiplication map, $\epsilon(1) = 1_A$ is the unit of the algebra, and $\eta(x) = B(1_A,x)$. Notice that the bilinear form $B$ can be also recovered from these data as $B(x,y) = \eta(xy)$.

In this section, we shall study a weak version of Frobenius algebras, as well as the problem of extending them to genuine Frobenius algebras.

\begin{defn}
An \emph{almost-Frobenius $k$-algebra} is a finite dimensional $k$-vector space $V$ endowed with a linear endomorphism $\tilde{T}: V \to V$, a linear map $\tilde{\epsilon}: k \to V$ and a linear functional $\tilde{\eta}: V \to k$.
\end{defn}

The rationale behind almost-Frobenius algebras is that they are the algebraic counterpart of finite-dimensional almost-TQFTs. 

\begin{prop}\label{prop:almost-frob-almost-tqft}
The category of (commutative) almost-Frobenius $k$-algebras is equivalent to the category of $(1+1)$ almost-TQFTs with values in $\Vecto{k}$.
\begin{proof}

Given an almost-Frobenius algebra $(V, \tilde{T}, \tilde{\epsilon}, \tilde{\eta})$, with the notation of Figure \ref{fig:decomp-sigmag}, let us define the $(1+1)$-dimensional almost-TQFT $Z$ given by $Z(S^1) = V$, $Z(D) = \tilde{\epsilon}$, $Z(D^\dag) = \tilde{\eta}$ and $Z(L) = \tilde{T}$. This automatically defines $Z$ on any morphism of $\Tubo{1+1}$, and thus of $\Tub{1+1}$ by monoidality.

Indeed, let $W: S^1 \to S^1 \in \Hom_{\Tubo{1+1}}(S^1, S^1)$ be a bordism with non-empty connected in and out boundaries. Then by gluing two discs to $W$ we get a closed orientable surface, which is diffeomorphic to $\Sigma_g$ for some $g \geq 0$. Hence, $W$ is diffeomorphic to $L^g: S^1 \to S^1$ and thus we set $Z(W) = \tilde{T}^g$. This representation of $W$ as a holed surface is unique, so the assignment is well-defined. The case in which $W$ has an empty boundary can be treated similarly.
%Furthermore, the process can be reversed, showing that this assignment is actually an equivalence. 
\end{proof}
\end{prop}

\subsection{Monoidal quantizability of almost-Frobenius algebras}\label{sec:monoidal-quant-almost-Frob}

Given an almost-Frobenius algebra $(V, \tilde{T}, \tilde{\epsilon}, \tilde{\eta})$, let us denote $v_0 = \tilde{\epsilon}(1) \in V$ and consider the collection of vectors $v_g = \tilde{T}^gv_0$ for $g \geq 0$. The vector subspace $W = \langle v_g\rangle_{g=0}^\infty$ of $V$ generated by these vectors will be called the \emph{core subspace} of $\tilde{T}$. We will say that the almost-Frobenius algebra is \emph{wide} if $W = V$. Observe that, since $V$ is finite dimensional, the subspace $W$ is generated by finitely many of these vectors $v_0, \ldots, v_{g_0}$. Notice that we can always get a wide algebra just by restricting $\tilde{T}$ and $\tilde{\eta}$ to $W$.

The following result shows that this information actually determines any potential Frobenius algebra extending an almost-Frobenius algebra.

\begin{figure}[h!]
\includegraphics[width=\textwidth]{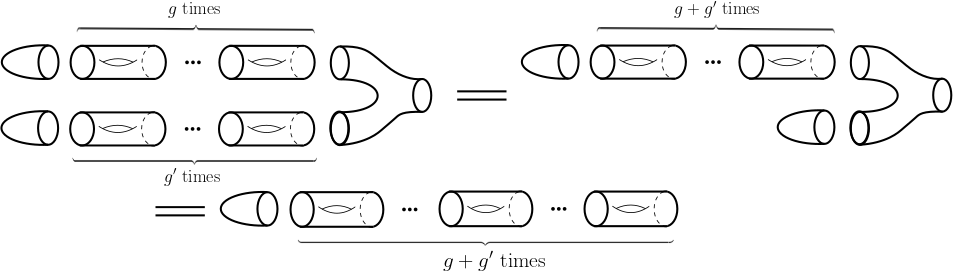}
\caption{The bordism $(T^gv_0) \cdot (T^{g'}v_0)$ is equivalent to $T^{g + g'}v_0$}
\label{fig:sum-tubes}
\centering
\end{figure}

\begin{prop}\label{prop:extension-almost-frob}
Let $(V, \tilde{T}, \tilde{\epsilon}, \tilde{\eta})$ be a wide almost-Frobenius algebra. There exists at most one Frobenius algebra structure on $V$ such that $T = \tilde{T}$, $\epsilon = \tilde{\epsilon}$ and $\eta = \tilde{\eta}$.

\begin{proof}
Suppose that there exists a Frobenius algebra on $V$ compatible with $\tilde{T}$. Then, since the Frobenius algebra satisfies $(T^gv_0) \cdot (T^{g'}v_0) = T^{g + g'}v_0$ (see Figure \ref{fig:sum-tubes}), then we must have that $v_g \cdot v_{g'} = v_{g+g'}$ for all $g,g' \geq 0$. Since the vectors $v_g$ span $V$, this characterizes uniquely the algebra structure and therefore also the bilinear form as $B(x,y) = \tilde{\eta}(xy)$.
\end{proof}
\end{prop}

\begin{rmk}
In the case that $(V, \tilde{T}, \tilde{\epsilon}, \tilde{\eta})$ is not wide, we still get a result by restricting to the core subspace $W \subseteq V$: There exists at most one Frobenius algebra structure on $W$ such that $T = \tilde{T}|_W$, $\epsilon = \tilde{\epsilon}|_W$ and $\eta = \tilde{\eta}|_W$. On the whole space, several Frobenius algebra structures may exist (see Example \ref{ex:c-frobenius}).
\end{rmk}

According to the previous result, we will say that a $\Vecto{k}$-valued $(1+1)$-dimensional almost-TQFT $Z$ is \emph{wide} if the associated almost-Frobenius algebra is wide. Notice that this means that the collection of vectors $v_g = Z(L^g \circ D)(1)$ for $g \geq 0$ spans $Z(S^1)$. Combining Propositions \ref{prop:almost-frob-almost-tqft} and \ref{prop:extension-almost-frob}, we get the following result.

\begin{cor}\label{cor:extension-almost-tqft}
Let $Z: \Tub{1+1} \to \Vecto{k}$ be a wide $(1+1)$-dimensional almost-TQFT. There exists at most one monoidal TQFT $\cZ: \Bord{1+1} \to \Vecto{k}$ extending $Z$, that is, such that the following diagram commutes.

\[\xymatrix{
\Bord{1+1} \ar@{--{>}}[r]^{\cZ\;\;\;\;} & \Vecto{k} \\
\Tub{1+1} \ar[u] \ar[ru]_{Z} &
}
\]
\end{cor}

Proposition \ref{prop:extension-almost-frob} and Corollary \ref{cor:extension-almost-tqft} not only are uniqueness results. They do provide an effective procedure for checking whether the almost-TQFT $Z$ can be promoted to a monoidal TQFT, since the (unique) algebra structure can be easily computed from the core submodule. 

\begin{thm}\label{thm:quantization-almostTQFT}
A $\Vecto{k}$-valued $(1+1)$-dimensional wide almost-TQFT $Z$ can be extended to a monoidal TQFT if and only if the following conditions hold:
\begin{enumerate}
	\item\label{cond:1} The bilinear form $B(v_g, v_{g'}) = Z(D^\dag)(v_{g+g'})$ is non-degenerate.
	\item\label{cond:2} Let $b_1, \ldots, b_N$ be any orthogonal basis of $Z(S^1)$ with respect to $B$. Then, for all $1 \leq j \leq N$, we have
	$$
		Z(L)(b_j) = \sum_{i=1}^N \frac{b_jb_i^2}{B(b_i, b_i)} .
	$$
	Here, for $b_i = \sum_g a^g_i\, v_{g}$ and $b_j = \sum_{g'} a^{g'}_j\, v_{g'}$, we set $b_i b_j = \sum_{g,g'} a_i^g a^{g'}_j\, v_{g+g'}$.
\end{enumerate}

\begin{proof}
The necessity of (\ref{cond:1}) is clear from the definition of a Frobenius algebra. For (\ref{cond:2}), observe that the counit of the Frobenius algebra is the element $\delta = \sum_i \frac{1}{B(b_i, b_i)} b_i \otimes b_i \in Z(S^1) \otimes_k Z(S^1)$, where $b_1, \ldots, b_N$ is any orthogonal basis of $Z(S^1)$. Hence, the comultiplication map is given by (c.f.\ Figure \ref{fig:comultiplication})
$$
	w(x) = (m \otimes \Id_{S^1})(\Id_{S^1} \otimes \delta)(x) = \sum_{i=1}^N \frac{1}{B(b_i, b_i)} (xb_i) \otimes b_i. 
$$
Therefore, the map $Z(L)$ is given by
$$
	Z(L)(x) = (m \circ w)(x) = \sum_{i=1}^N \frac{1}{B(b_i, b_i)} m\left(xb_i, b_i\right) = \sum_{i=1}^N \frac{1}{B(b_i, b_i)} xb_i^2.
$$

\begin{figure}[h!]
\includegraphics[width=7.5cm]{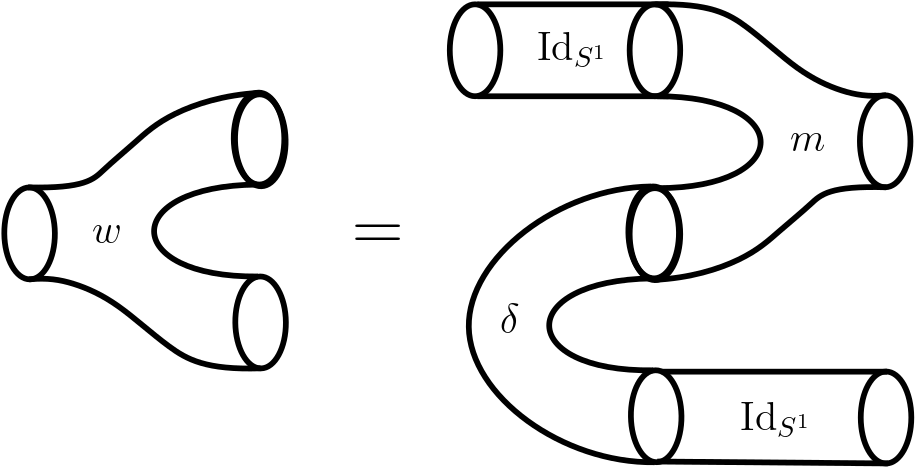}
\caption{The comultiplication map in terms of the elbow and the multiplication bordisms}
\label{fig:comultiplication}
\centering
\end{figure}

On the contrary, suppose that $Z$ is an almost-TQFT and let $(V, \tilde{T}, \tilde{\epsilon}, \tilde{\eta})$ be the almost-Frobenius algebra associated to $Z$ through Proposition \ref{prop:almost-frob-almost-tqft}. Since $V$ is wide, it is spanned by the vectors $v_0, \ldots, v_{N-1} \in V$, where $v_g = \tilde{T}^g\tilde{\epsilon}(1) \in V$ and, following Proposition \ref{prop:extension-almost-frob}, let us define a ring structure on $V$ by setting $v_{g} \cdot v_{g'} = v_{g + g'}$ for $g \geq 0$. It is clear that this multiplication is associative and $v_0$ is the unit of the ring. We also set as bilinear form $B(v_g, v_{g'}) = \tilde{\eta}(v_{g + g'})$.

Due to (\ref{cond:1}), this is a Frobenius algebra so it defines a monoidal TQFT $\cZ$. Furthermore, by construction $\cZ(D) = \tilde{\epsilon}$ and $\cZ(D^\dag) = \tilde{\eta}$, so in order to show that $\cZ$ extends $Z$ it only remains to prove that $\cZ(L) = Z(L)$. But, by the computation above, this is exactly condition (\ref{cond:2}).
\end{proof}
\end{thm}

\begin{rmk} Condition (\ref{cond:2}) fits perfectly with Proposition \ref{prop:necessity-quantization}. Indeed, if $Z$ satisfies (\ref{cond:2}), we have that
\begin{align*}
	Z(S^1 \times S^1)(1) &= Z(D^\dag) \circ Z(L) \circ Z(D)(1) = Z(D^\dag) \left(Z(L)(1) \right) \\
	&= Z(D^\dag) \left(\sum_{i=1}^N \frac{b_i^2}{B(b_i, b_i)} \right) = \sum_{i=1}^N \frac{Z(D^\dag)(b_i^2)}{B(b_i, b_i)}  = n.
\end{align*}
Where, in the last equality, we used that $Z(D^\dag)(b_i^2) = B(b_i,b_i)$.
\end{rmk}

\begin{rmk}
Recall that the monoidality condition of a monoidal TQFT forces it to take values in finite dimensional vector spaces. In this way, the condition that the almost-TQFT takes values on $\Vecto{k}$ is crucial for Theorem \ref{thm:quantization-almostTQFT}. Otherwise, no such extension can exist. 
\end{rmk}

Let us show how to check these conditions in practice in a couple of examples.

\begin{ex}
Consider the real almost-TQFT $Z: \Tub{1+1} \to \Vect{\RR}$ given by $Z(S^1) = \RR^2$ and
$$
	Z(L) = \begin{pmatrix} 0 & 1 \\ 1 & 2\end{pmatrix}, \quad Z(D)(1) = (1,0), \quad Z(D^\dag)(x,y) = x+y.
$$
Then, $v_0 = (1,0)$, $v_1=(0,1)$, $v_2=(1,2)$ so the core subspace coincides with the whole vector space $Z(S^1)$ and thus $Z$ is wide. In the basis $v_0, v_1$ the bilinear form is given by
$$
	B = \begin{pmatrix} 1 & 1 \\ 1 & 3\end{pmatrix},
$$
so applying the Gram-Schmidt orthogonalization procedure we get that $b_1 = v_0 = (1,0)$ and $b_2 = (1/\sqrt{2}, -1/\sqrt{2})$ is an orthonormal basis. Moreover, $b_2 \cdot b_2 = \frac12 (v_0 - v_1)^2 = (1,0)$ and, thus
$$
	\frac{1}{B(b_1, b_1)} b_1b_1^2 + \frac{1}{B(b_2, b_2)} b_1b_2^2 = (2,0) \neq Z(L)(b_1).
$$
Hence, $Z$ does not admit any extension to a monoidal TQFT. It is worth noticing that we actually knew that from Proposition \ref{prop:necessity-quantization}, since $Z(S^1 \times S^1)(1) = Z(\Sigma_1)(1) = Z(D^\dag)(v_1) = 1$ does not agree with the dimension of $Z(S^1)$.
\end{ex}

\begin{ex}\label{ex:c-frobenius}
In some cases, the argument in the proof of Proposition \ref{prop:extension-almost-frob} can be used even for non-wide TQFTs. Consider the real almost-TQFT $Z$ given by $Z(S^1) = \RR^2$ and
$$
	Z(L) = \begin{pmatrix} 2 & 0 \\ 0 & 2\end{pmatrix}, \quad Z(D)(1) = e_1 = (1,0), \quad Z(D^\dag)(x,y) = x.
$$
This is not a wide TQFT since the core subspace is $W = \langle (1,0) \rangle$. However, since $Z(D)(1) = e_1$, for any Frobenius algebra on $\RR^2$ extending this almost-Frobenius TQFT, the vector $e_1$ must be the unit of the ring. In particular, for $e_2 = (0,1)$ we have $e_1 \cdot e_2 = e_2$. Hence, any Frobenius algebra must satisfy
$$
	e_1 \cdot e_1 = e_1, \qquad e_1 \cdot e_2 = e_2, \qquad e_2 \cdot e_2 = a e_1 + b e_2,
$$
for certain $a, b \in \RR$. Hence, since $Z(D^\dag)(e_1) = 1$ and $Z(D^\dag)(e_2) = 0$, in this basis the bilinear form is
$$
	B = \begin{pmatrix} 1 & 0 \\ 0 & a\end{pmatrix}.
$$
Thus, we must have $a \neq 0$ and $e_1, e_2$ is an orthogonal basis so that the counit is $\delta = e_1 \otimes e_1 + \frac{1}{a} e_2 \otimes e_2$. Hence, we must have
$$
	2e_1 = Z(L)(e_1) = e_1^3 + \frac{1}{a} e_1e_2^2 = 2e_1 + \frac{b}{a}e_2, \quad 2e_2 = Z(L)(e_2) = e_2e_1^2 + \frac{1}{a} e_2^3 = be_1 + \left(2+ \frac{b^2}{a}\right)e_2.
$$
This forces $b = 0$. Hence, there exist infinitely many monoidal TQFTs extending the almost-TQFT $Z$, one for each value of $a \neq 0$. Up to isomorphism, this gives two possible Frobenius algebras: $\RR^2$ with $e_2 \cdot e_2 = 1$ (corresponding to $a=1$) and $\CC$ (corresponding to $a=-1$).
\end{ex}

\subsection{Almost-TQFTs as dynamical systems}\label{sec:dynamical-system}

In Section \ref{sec:monoidal-quant-almost-Frob} we have discussed what are the conditions for an almost-Frobenius algebra (equivalently, a $(1+1)$-dimensional almost-TQFT) to be promoted to a Frobenius algebra (i.e.\ a monoidal $(1+1)$-dimensional TQFT). However, in order to apply the previous results to quantize a diffeomorphism invariant of closed orientable surfaces, we still need to construct an initial almost-TQFT. To address this issue, the aim of this section is precisely to provide a method for constructing a wide almost-TQFT for an invariant (when possible). For this purpose, we shall adopt a different viewpoint by re-interpreting an almost-TQFT as a dynamical system.

Throughout this section, all the almost-TQFTs will be assumed to take values in the category $\Vecto{k}$ of finite dimensional $k$-vector spaces. Given categories $\cC$ and $\cD$, we will denote by $\textup{Fun}(\cC, \cD)$ the category of functors between $\cC$ and $\cD$, with natural transformations as morphisms. Moreover, if $\cC$ and $\cD$ are endowed with a monoidal structure, then $\textup{Fun}^{\otimes}(\cC, \cD)$ will denote the category of monoidal functors between $\cC$ and $\cD$. Finally, recall that $\Tubo{n+1}$ denotes the subcategory of $\Bord{n+1}$ of strict $(n+1)$-dimensional tubes and $\Tub{n+1}$ is its monoidal closure.

\begin{lem}\label{lem:dynam}
There exists an equivalence of categories
$$
	\textup{Fun}^{\otimes}(\Tub{n+1}, \Vecto{k}) \cong \textup{Fun}(\Tubo{n+1}, \Vecto{k})
$$
between $(n+1)$-dimensional almost-TQFTs and linear representations of the category $\Tubo{n+1}$.
\begin{proof}
The proof is immediate. Any almost-TQFT induces a representation of $\Tubo{n+1}$ by restriction. On the contrary, since $\Tub{n+1}$ is the closure under disjoint union of $\Tubo{n+1}$, the TQFT is obtained from a functor $\Tubo{n+1} \to \Vecto{k}$ by imposing monoidality.
\end{proof}
\end{lem}

\begin{cor}\label{cor:dynam}
Consider a presentation $\cT_{n+1}$ of $\Tubo{n+1}$. Then, there is an equivalence
$$
	\textup{Fun}^{\otimes}(\Tub{n+1}, \Vecto{k}) \cong \textup{Fun}^{\otimes}(\cT_{n+1}, \Vecto{k})
$$
between almost-TQFTs and linear representations of $\cT_{n+1}$. 
\end{cor}

Despite its simplicity, Corollary \ref{cor:dynam} leads to a natural interpretation of $\Vecto{k}$-valued almost-TQFTs in terms of linear discrete dynamical systems. A linear discrete dynamical system is a collection of finite dimensional $k$-vector spaces $\{V_\alpha\}_{\alpha \in \Lambda}$ (maybe with the indexing set $\Lambda$ infinite) and some linear maps $f_\beta: V_\alpha \to V_{\alpha'}$ between them. In this manner, given a functor $Z: \cT_{n+1} \to \Vecto{k}$, coming from an almost-TQFT, we define the dynamical system with vector spaces $\{Z(M)\}_{M \in \cT_{n+1}}$, one for each object in the presentation $\cT_{n+1}$, and maps $Z(W): Z(M_1) \to Z(M_2)$ for each morphism $W: M_1 \to M_2$ of $\cT_{n+1}$.

Furthermore, two linear dynamical systems $(V_\alpha, f_\beta)$ and $(V_\alpha', f_\beta')$ are equivalent (conjugate in the language of ODEs) if there exist linear isomorphisms $g_\alpha: V_\alpha \to V_\alpha'$ such that $f_\beta' = g_\alpha f_\beta g_\alpha^{-1}$. But this is nothing but the definition of a natural equivalence of almost-TQFTs, so classifying almost-TQFTs up to isomorphism is precisely to study linear discrete dynamical systems up to conjugacy. In particular, these observations open the door to studying almost-TQFT through the classical techniques of linear dynamical systems.

\begin{ex}
We have a presentation $\cT_{1+1}$ of $\Tubo{1+1}$ with two objects (corresponding to $\emptyset$ and $S^1$), two arrows (corresponding to the caps) and one loop (corresponding to $L: S^1 \to S^1$), as in the following diagram.
\[
\begin{tikzcd}
	\bullet \ar[bend right=60,swap]{r} &  \ar[bend right=60,swap]{l} \bullet \ar[loop, out=55 % start at angle 123°
    ,in= 305% stop at angle 57°
    ,distance=2.5em]{}{} 
\end{tikzcd}
\]
In this way, a $(1+1)$-dimensional almost-TQFT is the same as a vector space $V$ equipped with an endomorphism $V \to V$ and two maps $k \to V$ and $V \to k$. This provides an alternative proof of Proposition \ref{prop:almost-frob-almost-tqft} since this is precisely the information of an almost-Frobenius algebra.
\end{ex}

\begin{thm}\label{thm:almost-quantizable}
A non-constant invariant $\chi$ of closed orientable surfaces with values in a field $k$ is almost-quantizable through a $\Vecto{k}$-valued almost-TQFT if and only if there exist $N \geq 1$, a partition $(m_1, \ldots, m_s)$ of $N$, coefficients $\lambda_1, \ldots, \lambda_s \in \bar{k}$ in the algebraic closure $\bar{k}$ of $k$, and values $a_{i,j} \in \bar{k}$ for $1 \leq i \leq s$ and $0 \leq j \leq m_i-1$, such that for all $g \geq \max(m_1, \ldots, m_s)$
\begin{equation}\label{eq:invariant}
	\chi(\Sigma_g) = \sum_{i=1}^s \sum_{j=0}^{m_i-1} a_{i,j}\begin{pmatrix}g \\ j\end{pmatrix}\lambda_i^{g-j}.
\end{equation}

\begin{proof}
This formula is just the general solution to a linear discrete dynamical system. To make it more precise, notice that, without loss of generality, we can suppose that $k$ is algebraically closed. Suppose that there exists an almost-Frobenius algebra $(V, \tilde{T}, \tilde{\epsilon}, \tilde{\eta})$ of dimension $N$ computing $\chi$. As in Section \ref{sec:monoidal-quant-almost-Frob}, denote $v_0 = \tilde{\epsilon}(1)$ and $v_g = \tilde{T}^g(v_0)$ so that $\chi(\Sigma_g) = \tilde{\eta}(v_g)$. Since $k$ is algebraically closed, there exists an automorphism $P$ of $V$ such that $P\tilde{T}P^{-1} = J$ is of Jordan type. Then, the closed formula (\ref{eq:invariant}) follows from taking $(m_i, \ldots, m_s)$ as the partition of $N$ defined by the Jordan blocks of $J$ with eigenvalues $\lambda_1, \ldots, \lambda_s$ and computing $\chi(\Sigma_g) = \tilde{\eta}(P^{-1}J^gPv_0)$.

On the contrary, if $\chi(\Sigma_g)$ satisfies (\ref{eq:invariant}) then we can find a linear form $\tilde{\eta}: k^N \to k$, a Jordan type matrix $J$ and a vector $v_0 \in k^N$ such that $\chi(\Sigma_g) = \tilde{\eta}(J^gv_0)$. Setting $V = k^N$, $\tilde{T} = J$ and $\tilde{\epsilon}(1) = v_0$, we get the desired almost-Frobenius algebra.
\end{proof}
\end{thm}

\begin{rmk}\label{rmk:almost-quantizable-diagonal}
If the endomorphism $\tilde{T}$ is diagonalizable, then formula (\ref{eq:invariant}) takes the simpler form
$$
	\chi(\Sigma_g) = a_1\lambda_1^g + a_2\lambda_2^g + \ldots + a_N\lambda_N^g.
$$
This happens, for instance, if the bilinear form $B(v_g, v_{g'}) = \tilde{\epsilon}(v_{g+g'})$ is non-degenerate, as in condition (\ref{cond:1}) of Theorem \ref{thm:quantization-almostTQFT}.
\end{rmk}

\begin{rmk}\label{rmk:iteration}
An almost-quantizable invariant $\chi$ of closed orientable surfaces can be understood as a higher order recurrence. Let $(V, \tilde{T}, \tilde{\epsilon}, \tilde{\eta})$ be an almost-Frobenius algebra almost-quantizing $\chi$ and let $v_g = \tilde{T}^g\tilde{\epsilon}(1) \in V$ for $g \geq 0$. If $\dim V = N$, then there exist coefficients $a_0, \ldots, a_{N-1} \in k$ such that
$$
	v_N = a_0v_0 + a_1v_1 + \ldots + a_{n-1} v_{N-1}.
$$
Furthermore, since $\tilde{T}^k v_i = v_{i+k}$, we have that for all $k \geq 0$
$$
	v_{N+k} = \tilde{T}^k (v_N) = \sum_{i=0}^{N-1} a_i \tilde{T}^k (v_i) = \sum_{i=0}^{N-1} a_i  v_{i+k}.
$$
In particular, since $\chi(\Sigma_g) = \tilde{\eta}(v_g)$ we have that
$$
	\chi(\Sigma_{N+k}) = \sum_{i=0}^{N-1} a_i \chi(\Sigma_{i+k}),
$$
which is a recurrence relation of order $N$ for the invariant.
\end{rmk}

\begin{rmk}\label{rmk:computation-tqft}
If the invariant $\chi(\Sigma_g)$ is almost-quantizable in $\Vecto{k}$ as in Theorem \ref{thm:almost-quantizable}, there is an easy way of computing an almost-Frobenius algebra almost-quantizing it. Given $g\geq 0$ and $N \geq 1$, let us form the vectors
$$
	v_{g, N} = \left(\chi(\Sigma_g), \chi(\Sigma_{g+1}), \ldots, \chi(\Sigma_{g+N -1 })\right) \in k^N.
$$
As explained in Remark \ref{rmk:iteration}, there exists a maximum integer $N \geq 0$ such that the vectors $v_{0, N}, v_{1, N}, \ldots, v_{N-1, N}$ of $k^N$ are linearly independent. Hence, they form a basis of $k^N$ so we can write
$$
	v_{N,N} = a_0 v_{0,N} + a_1 v_{1,N} + \ldots + a_{N-1} v_{N-1,N},
$$
for some constants $a_0, \ldots, a_{N-1} \in k$. The coincidence of notation is not accidental: these coefficients are precisely the recurrence coefficients of Remark \ref{rmk:iteration}. In this form, the almost-Frobenius algebra can be easily recovered. Indeed, $V = k^N$, $\tilde{\epsilon}(1) = v_{0,N}$ and $\tilde{\eta}: k^N \to k$ is the projection onto the first component. The linear map $\tilde{T}$ satisfies, in the basis $v_{0, N}, v_{1, N}, \ldots, v_{N-1, N}$,
$$
	\tilde{T}(v_{i,N}) = \tilde{T}(v_{i+1, N}) \textrm{ for } 0 \leq i \leq N-1, \qquad \tilde{T}(v_{N-1,N}) = a_0 v_{0,N} + a_1 v_{1,N} + \ldots + a_{N-1} v_{N-1,N}.
$$
Hence, only with the knowledge of the invariant on the surfaces of genus $\leq 2N-1$, we can recover the whole almost-Frobenius algebra, and thus compute it for arbitrary genus.
\end{rmk}

\section{TQFTs with basepoints and representation varieties}\label{sec:representation-varieties}

In this section, we shall apply the previous results to the case of representation and character varieties. 

\begin{defn}
Let $G$ be an algebraic group and let $M$ be a compact orientable manifold. The \emph{$G$-representation variety} of $M$ is the variety parametrizing representations of the fundamental group $\pi_1(M)$ of $M$ into $G$
$$
	\Rep{G}(M) = \Hom\left(\pi_1(M), G\right).
$$
The Geometric Invariant Theory (GIT) quotient of $\Rep{G}(M)$ under the adjoint action of $G$ is the so-called \emph{character variety}, denoted by
$\Char{G}(M) = \Rep{G}(M) \sslash G$.
\end{defn}

\begin{rmk}
The algebraic structure of $\Rep{G}(M)$ is induced as follows. Let us consider a presentation of $\pi_1(M)$ with finitely many generators $\pi_1(M) = \langle \gamma_1, \ldots, \gamma_s \,|\, R_\alpha(\gamma_1, \ldots, \gamma_s) = 1\rangle$. Hence, we have a natural identification of $\Rep{G}(M)$ with the algebraic set $
\{(g_1, \ldots, g_s) \in G^s \,|\, R_\alpha(g_1, \ldots, g_s) = 1\}$. Moreover, different presentations give rise to isomorphic algebraic sets, so this induces a unique algebraic structure on $\Rep{G}(M)$.
\end{rmk}

Let us consider an algebraic invariant of algebraic varieties, that is a map $\phi: \Var{k} \to R$ into some (commutative, unitary) ring $R$ such that $\phi(X) = \phi(Y)$ if $X$ is isomorphic to $Y$. We want to focus on those invariants that are subtle enough to capture the dimension of the algebraic varieties.

\begin{defn}\label{defn:invariant}
An algebraic invariant $\phi: \Var{k} \to R$ is said to be \emph{weakly dimension-aware} if there exists a group homomorphism $d: (R-\{0\}, \cdot) \to (\ZZ, +)$, called the \emph{dimension witness}, and a natural number $N \geq 1$ such that $d(n) = 0$ for all $n$ in the ring of integers of $R$ and $d(\phi(X)) = N\dim X$ for every smooth algebraic variety $X$. It is said to be \emph{dimension-aware} if this equality holds for any algebraic varieties, even if they are not smooth.
\end{defn}

\begin{ex}
The Poincar\'e polynomial $P_c(X) = \sum_k (-1)^k H_c(X) t^k \in \ZZ[t]$ of the compactly supported (sheaf) cohomology of an algebraic variety $X$ is a weakly dimension-aware invariant since, for $X$ smooth, the degree of $P_c(X)$ is twice the dimension of $X$.
\end{ex}

\begin{ex}
An important invariant associated to an algebraic variety is its virtual class in the Grothendieck ring $\K{\Var{k}}$ of algebraic varieties. Recall that $\K{\Var{k}}$ is the ring generated by isomorphism classes of algebraic varieties, denoted by $[X]$, modulo the cut-and-paste relations
$
	[X] = [Y] + [X - Y]
$, for any algebraic variety $X$ and any locally closed set $Y \subseteq X$. In this way, we can form the invariant $\phi: \Var{k} \to \K{\Var{k}}$, $\phi(X) = [X]$, called the \emph{virtual class} or the \emph{motive} of $X$. Notice that the Poincar\'e polynomial with compact support factorizes through the Grothendieck ring of algebraic varieties
\[\xymatrix{
\Var{k}^{sm} \ar[r]^{P_c} \ar[d]_{\phi} & \ZZ[t]\\
\K{\Var{k}} \ar@{--{>}}[ru]
}
\]
Here, $\Var{k}^{sm}$ denotes the subcategory of $\Var{k}$ of smooth algebraic varieties. In particular, since the virtual class of any algebraic variety can be written as the sum of smooth subvarieties \cite{bittner2004universal}, this implies that the virtual class of an algebraic variety is a dimension-aware invariant.
\end{ex}

\begin{ex}
In the case that the underlying field is $k = \CC$, we can consider a coarser invariant from the Grothendieck ring. The compactly-supported cohomology of a complex algebraic variety $X$ is naturally equipped with a mixed Hodge structure $(H^k_c(X), W_\bullet, F^\bullet)$, where $W_\bullet$ is an increasing filtration of $H^k_c(X)$, called the weight filtration, and $F^\bullet$ is a decreasing filtration of the complex-valued cohomology $H^k_c(X; \CC)$, called the Hodge filtration (see \cite{Peters-Steenbrink:2008}). Using these filtrations, we can define the compactly-supported Hodge numbers of $X$ as $h^{k;p,q}_c(X) = \dim \textrm{Gr}^F_p \left(\left(\textrm{Gr}_W^{p+q} H_c^k(X) \right)\right) \otimes_\QQ \CC$.

These numbers can be stacked into the so-called \emph{$E$-polynomial} of $X$, also known as the Deligne-Hodge polynomial of $X$
$$
	e(X) = \sum_{k,p,q} (-1)^k h_c^{k; p, q}(X) \, u^p v^q \in \ZZ[u,v].
$$
Thanks to the splitting properties of the compactly-supported cohomology, we get that the $E$-polynomial respects disjoint unions and multiplications, giving rise to a ring homomorphism
$$
	e: \K{\Var{\CC}} \to \ZZ[u,v]
$$
In the case that $h_c^{k;p,q}(X) = 0$ for $p \neq q$, the $E$-polynomial only depends on the product $uv$ and it is customary to write it in the new variable $q = uv$. If $X$ is a smooth irreducible variety of complex dimension $n$, Poincar\'e duality shows that $h_c^{2n; 0, 0}(X) = 1$ and $h_c^{2n; p, q}(X) = 0$ for $(p,q) \neq (0,0)$. Moreover $h_c^{k; p, q}(X) = 0$ for $p+q > k$, showing that the degree of the $E$-polynomial is twice the (complex) dimension of $X$ for smooth irreducible varieties and thus, by additivity, for all varieties. Hence, the $E$-polynomial is a dimension-aware invariant.
\end{ex}

Let $\phi$ be a dimension-aware invariant and $G$ an algebraic group, and let us form the invariant $\chi(M) = \phi(\Rep{G}(M))$. The following result shows that there is no hope of quantizing $\chi$ through a monoidal TQFT.

\begin{thm}\label{thm:non-quantizable-rep-var}
If $G$ is an algebraic group with $\dim G > 0$, then any dimension-aware invariant of the $G$-representation variety is not a strongly quantizable invariant.
\begin{proof}
Consider the manifold $M = S^{n} \times S^1$ (in the case $n=0$, we take $M = S^1$). For $n \neq 1$, we have $\pi_1(M) = \pi_1(S^1) = \ZZ$ and, thus, $\Rep{G}(M) = G$ so in particular $\dim \Rep{G}(M) > 0$. Hence, for any dimension-aware invariant $\phi$, we cannot have $\phi(\Rep{G}(M))$ lying in the ring of integers of $R$, as required by Corollary \ref{cor:non-quantizable-integers}.

In the case $n=1$, use that we have a diagonal embedding $G \hookrightarrow \Rep{G}(S^1 \times S^1)$, which also implies that $\dim \Rep{G}(S^1 \times S^1) > 0$, and we proceed as above.
\end{proof}
\end{thm}

\subsection{Quantization with basepoints}\label{sec:tqfts-with-basepoints}

As shown in Theorem \ref{thm:non-quantizable-rep-var}, the virtual classes of representation varieties are not strongly quantizable. However, a closer look evidences that this is an expectable fact, since the fundamental group is not a functor out of the category of topological spaces, but of pointed topological spaces. Hence some sort of extra data seems to be needed. For this purpose, we shall enlarge our category of bordisms by adding a finite collection of basepoints. This new category $\Bordp{n+1}$ is comprissed of the following information:
\begin{itemize}
	\item Objects: The objects of $\Bordp{n+1}$ are pairs $(M, A)$ where $M$ is an object of $\Bord{n+1}$ and $A \subseteq M$ is a finite set of basepoints with non-empty intersection with every connected component of $M$.
	\item Morphisms: A morphism $(M_1, A_1) \to (M_2, A_2)$ is a class of pairs $(W, A)$ where $W$ is a $(n+1)$-dimensional bordism between $M_1$ and $M_2$ as in $\Bord{n+1}$ and $A \subseteq W$ is a non-empty finite set such that $A \cap M_1 = A_1$ and $A \cap M_2 = A_2$. Two such pairs $(W, A), (W', A')$ are declared as equivalent if there exists a boundary-preserving diffeomorphism $F: W \to W'$ with $F(A)=A'$.
	\item Composition: Given two bordisms of the form $(W, A): (M_1, A_1) \to (M_2, A_2)$ and $(W', A'): (M_2, A_2) \to (M_3, A_3)$, the composed morphism $(W', A') \circ (W, A): (M_1, A_1) \to (M_3, A_3)$ is the class of the pair $(W \cup_{M_2} W', A \cup A')$. 
	\item Monoidality: $\Bordp{n+1}$ is endowed with its natural monoidal structure given by disjoint union (of both manifolds and basepoints).
\end{itemize}
For simplicity, the singleton set of basepoints will be denoted by $\star$. Thus $(M, \star)$ will be a manifold with a single marked point, which we loosely will denote by $M$. Notice that, to obtain a genuine category, we also need to add `thin' artificial bordisms $(M, A): (M,A) \to (M,A)$ for any object $(M,A)$ as identity arrows.

\begin{defn}
Let $\cC$ be a symmetric monoidal category. A \emph{monoidal (resp.\ lax monoidal) TQFT with basepoints} with values in $\cC$ is a monoidal (resp.\ lax monoidal) symmetric functor
$$
	\cZ: \Bordp{n+1} \to \cC.
$$
\end{defn}

\begin{rmk}\label{rmk:almost-tqft-from-basepoints}
In the $(1+1)$-dimensional case, we have a faithful functor $\kappa^0: \Tubo{1+1} \to \Bordp{1+1}$ from $2$-dimensional strict tubes without basepoints to $2$-dimensional bordisms with basepoints. On objects, it just sends $\kappa^0(S^1) = (S^1, \star)$. For a tube $\Sigma$, the functor assigns $\kappa^0(\Sigma) = (\Sigma, A_{\Sigma})$, being $A_{\Sigma}$ a set of cardinality
$$
	|A_{\Sigma}| = - \frac{\chi(\Sigma) + h^0(\partial\Sigma) - 4}{2},
$$
where $\chi$ denotes the Euler characteristic and $h^0(\partial\Sigma)$ is the number of components of the boundary of $\Sigma$. This set $A_\Sigma$ has one point in each connected components of $\partial\Sigma$ and the remaining ones lie in the interior of $\Sigma$. Notice that $\chi(\Sigma) + h^0(\partial\Sigma)$ is the Euler characteristic of the surface with $h^0(\partial\Sigma)$ discs attached to the boundaries of $\Sigma$, so the fact that $\kappa^0$ preserves composition follows from the usual formula for the Euler characteristic of the connected sum of surfaces. In particular, for $\Sigma_g$ the closed genus $g$ surface, the bordism $\kappa^0(\Sigma_g)$ has $g+1$ basepoints.

In this way, we get a monoidal functor $\kappa: \Tub{1+1} \to \Bordp{1+1}$. Therefore, this functor implies that any $(1+1)$-dimensional lax monoidal TQFT with basepoints $\cZ: \Bordp{1+1} \to \cC$ induces an almost-TQFT without basepoints $Z = \cZ \circ \kappa: \Tub{1+1} \to \cC$.
\end{rmk}

Observe that the question of quantizability for TQFTs with basepoints is subtler than for usual TQFTs since $\Bordp{n+1}$ contains more morphisms than $\Bord{n+1}$. In particular, we cannot directly apply Proposition \ref{prop:necessity-quantization} to give necessary conditions to the quantizability of algebraic invariants through TQFTs with basepoints.

To set this problem, in analogy with the unpointed case of Section \ref{sec:quantizability-problem}, a $\Lambda$-valued diffeomorphism invariant $\chi$ of manifolds with basepoints will be an assignment of an element $\chi(W, A) \in \Lambda$ to any pair of a closed orientable manifold $W$ and a non-empty finite subset $A \subseteq W$ such that, if $f: W \to W'$ is a diffeomorphism with $f(A) = A'$, then $\chi(W, A) = \chi(W', A')$. 

As we will see, the basepoints prevent us to have a duality property on the nose, as used in the proof of Proposition \ref{prop:necessity-quantization}. To get some control on the basepoints, we introduce the following special bordism.

\begin{defn}
Given a closed connected orientable manifold $M$, the \emph{twisting bordism} is the morphism
$$
	\Theta_M = (M \times [0,1], \{x_0, x_1\}): M \to M,
$$
where $x_0 \in M \times \{0\}$ and $x_1 \in M \times \{1\}$, as depicted in Figure \ref{fig:twist}.

\begin{figure}[h!]
\includegraphics[width=7cm]{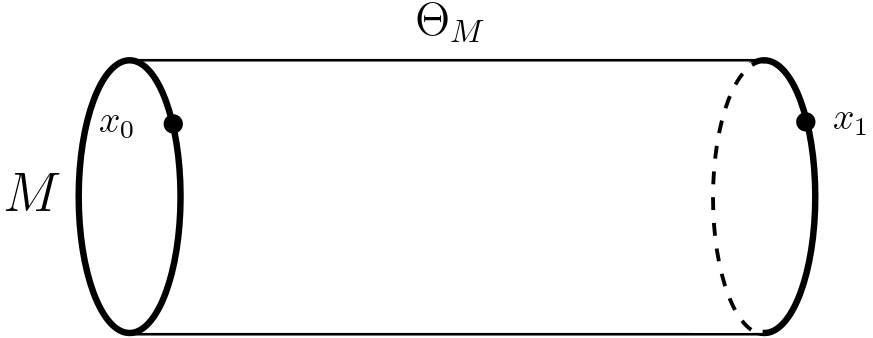}
\caption{The twisting bordism}
\label{fig:twist}
\centering
\end{figure}
\end{defn}

In this setting, keeping track of the basepoints in the proof of Proposition \ref{prop:necessity-quantization}, we get the following result in the case $\cC = \Mod{R}$.

\begin{prop}\label{prop:trace-twisting}
Consider a monoidal TQFT with basepoints
$$
	\cZ: \Bordp{n+1} \to \Mod{R}.
$$
Then, for any closed connected orientable $n$-dimensional manifold $M$, we have that the induced homomorphism $\cZ(M \times S^1, \{x_1, x_2\}): R \to R$ 
where $x_1, x_2$ are any two basepoints, is given by multiplication by the trace $\tr(\cZ(\Theta_M^2)|_V)$ of $\cZ(\Theta_M^2)$ restricted to a finitely generated invariant submodule $V \subseteq \cZ(M, \star)$.
\begin{proof}
Let $\delta: \emptyset \to (M, \star) \sqcup (M, \star)$ and $\mu: (M, \star) \sqcup (M, \star) \to \emptyset$ be the two `elbow' bordisms with two basepoints. As depicted in Figure \ref{fig:zorro-punctures}, in this setting we have that $(\mu \sqcup \Id_M) \circ (\Id_M \sqcup \delta) = (M \times [0,1], \{x, x', x''\}) = \Theta_M^2$.

\begin{figure}[h!]
\includegraphics[width=9cm]{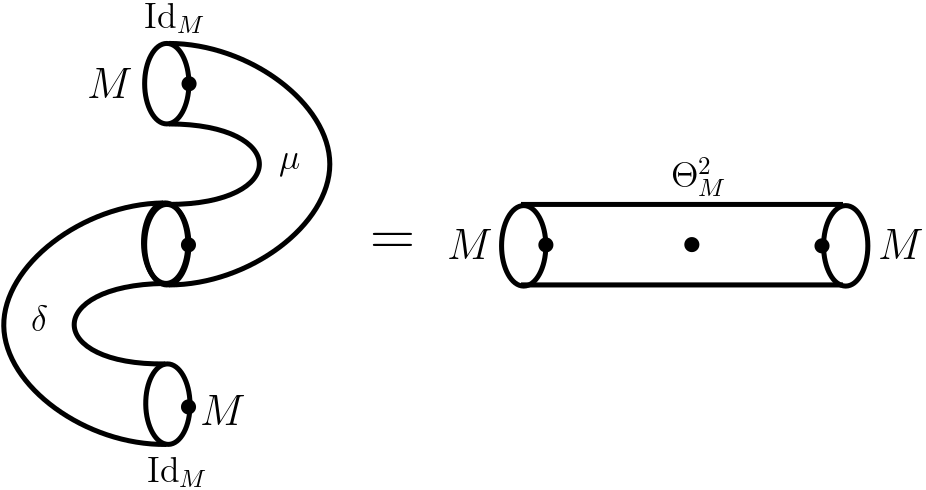}
\caption{The Zorro argument with punctures}
\label{fig:zorro-punctures}
\centering
\end{figure}

Set $\cZ(\delta)(1) = \sum_{i,j=1}^N a_{i,j} \,v_i \otimes v_j$ with $v_1, \ldots, v_N \in \cZ(M, \star)$.
Hence, under the TQFT, we have that
\begin{align}\label{eq:trace}
	\cZ(\Theta_M^2)(v_k) &=  (\cZ(\mu) \otimes \cZ(\Id_M)) \circ (\cZ(\Id_M) \otimes \cZ(\delta))(v_k) = \sum_{i,j=1}^N a_{i,j} \cZ(\mu)(v_k, v_i) v_j.
\end{align}
Therefore, we have that
\begin{align*}
	\cZ(M \times S^1)(1) & =  \cZ(\mu) \circ \cZ(\delta)(1) = \cZ(\mu) \left( \sum_{i,j=1}^N a_{i,j} v_i \otimes v_j\right) = \sum_{i,j=1}^N a_{i,j} \cZ(\mu)(v_i, v_j). 
\end{align*}
This is precisely the trace of $\cZ(\Theta_M^2)$ restricted to the submodule $V \subseteq \cZ(M, \star)$ generated by $v_1, \ldots, v_N$, as can be checked from (\ref{eq:trace}).
\end{proof}
\end{prop}

\begin{defn}
An $R$-valued invariant $\chi$ is said to be \emph{split} with factor $\theta \in R$ if, for any closed connected orientable manifold $W$ and any non-empty finite set $A \subseteq W$, we have
$$
	\chi(W, A \cup \star) = \theta \, \chi(W, A).
$$
Additionally, a TQFT $\cZ: \Bordp{n+1} \to \Mod{R}$ is said to be \emph{split} with factor $\theta \in R$ if for any bordism $(W,A)$ of $\Bordp{n+1}$ we have
$$
	\cZ(W, A \cup \star) = \theta \cZ(W, A).
$$
\end{defn}

\begin{rmk}
Since $\Theta_M^2 =  (M \times [0,1], \{x_0, x_1\} \cup \star)$, in particular $\cZ(\Theta_M^2) = \theta \cZ(\Theta_M)$.
\end{rmk}

\begin{rmk}
Using several times the definition of being split, for any split invariant $\chi$ with factor $\theta$ and any non-empty finite set $A \subseteq W$ we have that
$$
	\chi(W, A) = \theta^{|A|-1} \chi(W, \star).
$$
Analogously, if $\cZ$ is a split TQFT then
$$
	\cZ(W, A) = \theta^{|A|-|A_0|} \cZ(W, A_0),
$$
where $A_0 \subseteq W$ is a finite set with a single point on each connected component of $\partial W$. Notice that, by the definition of a morphism in $\Bordp{n+1}$, no more points can be removed.
\end{rmk}

\begin{rmk}\label{rmk:wide-implies-split}
Not necessarily a monoidal TQFT $\cZ$ computing a split invariant is split itself. However, the situation is different if $\cZ$ is wide, i.e.\ for any object $M$, the $R$-module $\cZ(M)$ is generated by the images $\cZ(B)(1)$ of the possible ``fillings'' $B: \emptyset \to M$. In that case, for any bordism $W: M_1 \to M_2$ we have that the non-degenerate bilinear map $(v_1, v_2) \mapsto \cZ(\mu_{M_2})(\cZ(W \circ \Theta_{M_1})(v_1), v_2)$ equals the bilinear map $(v_1, v_2) \mapsto \theta\cZ(\mu_{M_2})(\cZ(W)(v_1), v_2)$ and thus $\cZ(W \circ \Theta_{M_1}) = \theta\cZ(W)$, so $\cZ$ is split.
\end{rmk}

\begin{lem}\label{lem:eigen-twist}
If $\cZ$ is a split TQFT and $R$ is an integral domain, then the eigenvalues of $\cZ(\Theta_M)$ are either $\theta$ or $0$.

\begin{proof}
Since $\cZ(\Theta_M)^2 - \theta \cZ(\Theta_M) = 0$, then by the Cayley-Hamilton theorem we have that the roots of $p(x) = x^2-\theta x = x(x-\theta)$ are the possible eigenvalues of $\cZ(\Theta_M)$ in the ring of fractions of $R$.
\end{proof}
\end{lem}

\begin{cor}\label{cor:cond-non-quantizable-basepoints}
Let $\cZ: \Bordp{n+1} \to \Mod{R}$ be a split monoidal TQFT with basepoints taking values in $R$-modules, where $R$ is an integral domain. Then, for any closed connected orientable $n$-dimensional manifold $M$, we have that the induced homomorphism
$$
	\cZ(M \times S^1, \{x_1, x_2\}): R \to R,
$$ 
where $x_1, x_2$ are any two basepoints, is given by multiplication by $m\theta^2$ for some integer $m \geq 0$.

\begin{proof}
By Lemma \ref{lem:eigen-twist}, the eigenvalues of $\cZ(\Theta_M^2)$ are either $\theta^2$ or $0$, so $\tr(\cZ(\Theta_M^2)|_V) = m\theta^2$, where $m$ is the number of non-vanishing eigenvalues. Now the result follows from Proposition \ref{prop:trace-twisting}.
\end{proof}
\end{cor}

\begin{cor}
Let $\chi$ be an split invariant of $(n+1)$-dimensional closed orientable manifolds with basepoints with values in an integral domain $R$ and factor $\theta\in R$. If there exists an $n$-dimensional closed orientable manifold $M$ such that $\chi(M \times S^1, \star) \not\in \theta\NN$, then $\chi$ is not strongly quantizable through a split monoidal TQFT. 
\end{cor}

\subsection{Split TQFTs for representation varieties}\label{sec:split-tqfts-rep}

Let us come back to the problem of quantizing representation varieties. Virtual classes of representation varieties can be extended to work with basepoints. Given a compact manifold $M$ and a finite set $A \subseteq M$, let $\Pi(M, A)$ be the fundamental groupoid of $M$ with basepoints in $A$, that is, the groupoid of homotopy classes of paths in $M$ between points of $A$. Then, given an algebraic group $G$, the representation variety of $(M, A)$ is the set of groupoid homomorphisms $\rho: \Pi(M, A) \to G$, that is,
$$
	\Rep{G}(M, A) = \Hom_{\Grpd}\left(\Pi(M, A), G\right).
$$

Decompose $M$ into connected components $M = M_1 \sqcup \ldots \sqcup M_s$ and let as choose a set $A_0 = \{x_1, \ldots, x_s\} \subseteq A$ with $x_i \in M_i$ for all $i$. Then, by choosing paths between the remaining points of $M_i \cap A$ and $x_i$, we have a natural identification
\begin{equation}\label{eq:product-rep-var}
	\Rep{G}(M, A) = \prod_{i=1}^s \left(\Hom(\pi_1(M_i), G) \times G^{|M_i \cap A| - 1} \right) = G^{|A|-s} \times \prod_{i=1}^s \Rep{G}(M_i).
\end{equation}
Hence, $\Rep{G}(M, A)$ is an algebraic variety. Moreover, we can consider the character variety of $(M,A)$ as the GIT quotient of $\Rep{G}(M,A)$ under the adjoint action of $G$, $\Char{G}(M,A) = \Rep{G}(M,A) \sslash G$.

Now, let us consider an invariant $\phi: \Var{k} \to R$ of algebraic varieties with values in a ring $R$, as in Definition \ref{defn:invariant}. We will say that $\phi$ is \emph{multiplicative} if $\phi(X \times Y) = \phi(X) \phi(Y)$ for all algebraic varieties $X$ and $Y$.

\begin{lem}\label{lem:multiplicative-split}
If $\phi$ is a multiplicative invariant with values in a ring $R$, then the $\chi(M, A) = \phi(\Rep{G}(M,A))$ is a split invariant with factor $\phi(G) \in R$.

\begin{proof}
It is a straightforward computation. Let $M$ be a closed connected manifold and let $A \subseteq W$ be a finite non-empty set. By (\ref{eq:product-rep-var}) we have that
\begin{align*}
	\Rep{G}(W, A \cup \star) &= \Hom_{\Grpd}(\Pi(W, A \cup \star), G) \\
	&= \Hom_{\Grpd}(\Pi(W, A), G) \times G = \Rep{G}(W, A) \times G.
\end{align*}
In particular, $\phi(\Rep{G}(W, A \cup \star)) = \phi(G) \phi(\Rep{G}(W, A))$, as we wanted to prove.
\end{proof}
\end{lem}

\begin{prop}\label{prop:dim-rep-var-torus}
Suppose that $k$ is an algebraically closed field. If $G$ is an algebraic group over $k$, then $\dim \Rep{G}(S^1 \times S^1) > \dim G$.
\begin{proof}
Notice that we have
$$
	\Rep{G}(S^1 \times S^1) = \{(g,h) \in G^2 \,|\, [g,h]= 1\} \cong  \{(g,h) \in G^2 \,|\, hgh^{-1}= g\}.
$$	
In this way, the projection onto the first component $\Rep{G}(S^1 \times S^1) \to G$ has as fiber the centralizer $C_{g}(G)$ of $g$ inside $G$. Hence, if $\dim C_{g}(G) > 0$ for all $g$ in a dense open set of $G$, the result follows.

We actually have that $\dim C_{g}(G) > 0$ for any $g \in G$\footnote[2]{This proof was communicated to us by F. Knop.}. Indeed, $\dim C_{g}(G) \geq \dim C_{g}(B)$, where $B$ is a Borel subgroup containing $g$. Now, the conjugacy class of $g$ in $B$ is contained in the coset $g^{-1}(B,B)$, with $(B,B)$ the commutator subgroup of $B$. But $\dim (B, B) < \dim B$ since $B$ is solvable (and non-trivial), so the dimension of the conjugacy class is also strictly smaller than $\dim B$. Thus, $\dim C_{g}(B) > 0$, as we wanted to prove.
\end{proof}
\end{prop}

\begin{thm}\label{thm:non-quantizable-basepoints-rep-var}
If $G$ is an algebraic group with $\dim G > 0$, then any multiplicative dimension-aware invariant with values in an integral domain of the $G$-representation variety of manifolds of dimension $\geq 2$ is not a strongly quantizable invariant through a split monoidal TQFT with basepoints.

\begin{proof}
Let us first consider the case of surfaces. Let $\phi$ be the dimension-aware invariant and suppose that there exists a monoidal TQFT with basepoints $\cZ: \Bordp{1+1} \to \Mod{R}$ quantizing the $\phi$-invariant of the $G$-representation variety. By Corollary \ref{cor:cond-non-quantizable-basepoints} and Lemma \ref{lem:multiplicative-split}, we have that
$$
	m\phi(G^2) = \cZ(S^1 \times S^1) = \phi(\Rep{G}(S^1 \times S^1, \{x_1,x_2\})) = \phi(G \times \Rep{G}(S^1 \times S^1)),
$$
for some natural number $m \geq 0$. Hence, since $\phi$ is dimension-aware and multiplicative, we get
\begin{align*}
N(\dim G + \dim \Rep{G}(S^1 \times S^1)) &= d(\phi(G \times \Rep{G}(S^1 \times S^1))) = d(m\phi(G^2)) \\
&= d(m) + d(\phi(G^2)) = d(\phi(G^2)) = 2 N\dim G,
\end{align*}
contradicting Proposition \ref{prop:dim-rep-var-torus}.

For dimension $n +1 > 3$, the same argument as above works verbatim with the manifold $S^{n-1} \times S^1 \times S^1$. In the case $n = 2$, notice that $ \Rep{G}(S^1 \times S^1)$ is a subvariety of $ \Rep{G}(S^1 \times S^1 \times S^1)$. Hence, the previous argument implies that $2 \dim G = \dim G + \dim \Rep{G}(S^1 \times S^1 \times S^1) \geq \dim G + \dim \Rep{G}(S^1 \times S^1)$, which is again a contradiction.
\end{proof}
\end{thm}

\begin{rmk}
The hypothesis of being a split TQFT in Theorem \ref{thm:non-quantizable-basepoints-rep-var} can be substituted for being a wide TQFT, as pointed out in Remark \ref{rmk:wide-implies-split}.  
\end{rmk}

\begin{cor}\label{cor:non-quantizable-basepoints-rep-var}
If $G$ is an algebraic group with $\dim G > 0$, then the virtual class of the $G$-representation variety of manifolds of dimension $\geq 2$  is not a strongly quantizable invariant through a split monoidal TQFT with basepoints.
\end{cor}

\begin{proof}
Let $R = \K{\Var{k}}$ be the Grothendieck ring of algebraic varieties. Notice that $R$ is not a domain, so we cannot directly apply Theorem \ref{thm:non-quantizable-basepoints-rep-var}. However, the Poincar\'e polynomial with compact support $P_c: \K{\Var{k}} \to \ZZ[t]$ is a ring homomorphism into an integral domain. Now, let $f = \cZ(\Theta_M): \cZ(M) \to \cZ(M)$ be the image of the twisting bordism of $M$ under a monoidal TQFT. By fixing a set of generators, we have lifting maps
\[\xymatrix{
\ZZ[t]^s  \ar@{--{>}}[r]^{f''} & \ZZ[t]^s \\ 
R^s \ar[u]^{P_c^s} \ar[d] \ar@{--{>}}[r]^{f'} & R^s \ar[d] \ar[u]_{P_c^s} \\
\cZ(M) \ar[r]^{f} & \cZ(M) \\
}
\]
Since $\ZZ[t]$ is an integral domain, the trace of $f''$ is $mP_c(G)^2$ for some integer $m \geq 0$. Now, since the witness of the dimension of the virtual class is the degree $d([X]) = \deg P_c([X])$, the result follows by using the same argument as in Theorem \ref{thm:non-quantizable-basepoints-rep-var}.
\end{proof}

\begin{rmk}
The previous argument can be applied to other related situations, such as the virtual classes of $G$-character stacks as studied in \cite{gonzalez2022virtual}, to show that they are not strongly quantizable either.
\end{rmk}

\subsection{Lax monoidal TQFTs for representation varieties}\label{sec:lax-monoidal-rep-var}

We have seen in Theorem \ref{thm:non-quantizable-basepoints-rep-var} that the virtual class of $G$-representation varieties cannot be strongly quantized through a monoidal TQFT, even if we consider basepoints. However, it turns out that it can be lax-quantized through a lax monoidal TQFT with basepoints.

\begin{thm}{\cite[Theorem 4.9]{GP-2018}}\label{thm:lax-monoidal-rep-var}
Let $G$ be an algebraic group and $n \geq 0$. There exists a split lax monoidal TQFT with basepoints
$$
	\cZ_G: \Bordp{n+1} \to \Mod{\K\Var{k}}
$$
such that $\cZ_G(W, \star)(1) = [\Rep{G}(W)] \in \K{\Var{k}}$ for any closed connected orientable $(n+1)$-dimensional manifold $W$. 
\end{thm}

For the sake of completeness, let us sketch briefly the construction of this TQFT. For further details, please refer to \cite{GPLM-2017,GP-2018}. Given an algebraic variety $X$ over $k$, let us denote by $\Varrel{X}$ the category of $X$-varieties, that is, the category whose objects are regular morphisms $f: Z \to X$ and whose morphisms between $f: Z \to X$ and $f': Z' \to X$ are intertwining regular maps $g: Z \to Z'$, i.e. satisfying $f' \circ g = f$. Given a regular map $h: X \to X'$ we have two induced functors, $h_!: \Varrel{X} \to \Varrel{X'}$ (covariant) and $h^*: \Varrel{X'} \to \Varrel{X}$ (contravariant), given by post-composition and pullback, respectively.

Passing to the Grothendieck ring, we get $\K{\Var{k}}$-algebras $\K{\Varrel{X}}$, and each regular morphism $h: X \to X'$ induces a $\K{\Var{k}}$-module homomorphism $h_!: \K{\Varrel{X}} \to \K{\Varrel{X'}}$ and an algebra homorphism $h^*: \K{\Varrel{X'}} \to \K{\Varrel{X}}$. With these notions, the functor $\cZ_G$ is given as follows.
\begin{itemize}
	\item On an object $(M, A) \in \Bordp{n+1}$, the functor assigns the $\K{\Var{k}}$-algebra $\cZ_G(M,A) = \K{\Varrel{\Rep{G}(M,A)}}$, seen as a  $\K{\Var{k}}$-module.
	\item For a bordism $(W, A): (M_1, A_1) \to (M_2, A_2)$, let us denote by $j_1: \Rep{G}(W,A) \to \Rep{G}(M_1,A_1)$ and $j_2: \Rep{G}(W,A) \to \Rep{G}(M_2,A_2)$ the restriction maps. Then, $\cZ(W, A) = (j_2)_! \circ (j_1)^*:  \K{\Varrel{\Rep{G}(M_1,A_1)}} \to  \K{\Varrel{\Rep{G}(M_2,A_2)}}$.
\end{itemize}
Using the Seifert-van Kampen theorem for fundamental groupoids and the properties of base change, it can be proven that $\cZ_G$ is actually a functor. Moreover, the external product $\boxtimes:  \K{\Varrel{\Rep{G}(M_1,A_1)}} \otimes  \K{\Varrel{\Rep{G}(M_2,A_2)}} \to  \K{\Varrel{\Rep{G}(M_1 \sqcup M_2, A_1 \sqcup A_2)}}$ endows $\cZ_G$ with the structure of a lax monoidal functor, giving rise to a lax monoidal TQFT.

\begin{rmk}
It is worth mentioning that, even though the virtual class of representation varieties of orientable surfaces are lax-quantizable, their dimensions are not due to their growth rate. By Theorem \ref{thm:almost-quantizable}, if the dimension of the representation variety was a quantizable invariant we would have that $\dim \Rep{G}(\Sigma_g)$ is given by Equation (\ref{eq:invariant}) for some coefficients $a_{i,j}, \lambda_i \in \overline{\QQ}$. This sequence is either constant or growths exponentially for $g \to \infty$. But $\Rep{G}(\Sigma_g) \subseteq G^{2g}$ so $\dim \Rep{G}(\Sigma_g) \leq 2g \dim G$. On the other hand $G^{g} \hookrightarrow \Rep{G}(\Sigma_g)$ as $(h_1, \ldots, h_g) \mapsto (h_1, 1, \ldots, h_g, 1)$ so $\dim \Rep{G}(\Sigma_g) \geq g \dim G$. This excludes both the constant and the exponential growth scenarios.
\end{rmk}

As shown in Theorem \ref{thm:non-quantizable-basepoints-rep-var}, if $\dim G \geq 1$, then the virtual classes of the $G$-representation varieties are not quantizable. Nevertheless, if $\dim G = 0$, then the necessary condition of Corollary \ref{cor:cond-non-quantizable-basepoints} may still hold, since $[\Rep{G}(S^1 \times S^1)] $ may be an integer, its number of points. In fact, in this case, the lax monoidal TQFT of Theorem \ref{thm:lax-monoidal-rep-var} is actually monoidal.

\begin{prop}\label{prop:TQFT-dim0-monoidal}
If $\dim G = 0$, then the TQFT of Theorem \ref{thm:lax-monoidal-rep-var} is monoidal.

\begin{proof}
Let $X$ be variety of dimension $0$. Given $x \in X$, let $1_x: \star \to X$ denote the inclusion map of $x$ into $X$. In this situation, we have a natural identification
$$
	\K{\Varrel{k}}^{\oplus |X|} \cong \K{\Varrel{X}}, \qquad (\alpha_x)_{x\in X} \mapsto \sum_{x\in X} \alpha_x 1_x,
$$
where $|X|$ denotes the number of closed points of $X$. 
In particular, using this isomorphism we get an identification $\K{\Varrel{(X \times Y)}} \cong \K{\Varrel{k}}^{\oplus |X||Y|} \cong \K{\Var{k}}^{\oplus |X|} \otimes_{\K{\Varrel{k}}} \K{\Varrel{k}}^{\oplus |Y|} \cong \K{\Varrel{X}} \otimes_{\K{\Var{k}}}  \K{\Varrel{Y}}$.

Now, if $\cZ_G: \Bordp{n+1} \to \Mod{\K\Var{k}}$ is the TQFT of Theorem \ref{thm:lax-monoidal-rep-var}, for any $n$-dimensional closed orientable manifolds $M$ and $M'$, we have $\cZ_G(M \sqcup M') = \K{\Varrel{(\Rep{G}(M) \times \Rep{G}(M'))}}$, which by the previous observation is isomorphic to $\K{\Varrel{\Rep{G}(M)}} \otimes_{\K{\Var{k}}} \K{\Varrel{\Rep{G}(M')}} = \cZ_G(M) \otimes_{\K{\Var{k}}} \cZ_G(M')$. An analogous computation with disjoint union of bordisms shows that they preserve this isomorphism.
\end{proof}
\end{prop}

\begin{rmk}
If $G$ is an algebraic group over $\FF_q$ (no necessarily of dimension $0$), then we can adapt the TQFT of Theorem \ref{thm:lax-monoidal-rep-var} to get a monoidal TQFT computing the number of points of the $G$-representation variety. For that purpose, we consider the forgetful functor of closed points $\Var{\FF_q} \to \Top$, $X \mapsto X(\FF_q)$, which gives rise to a ring homomorphism $\K{\Var{\FF_q}} \to \K{\Top}$ and its relative versions $\K{\Varrel{Z}} \to \K{\Top/Z(\FF_q)}$. These maps induce a natural transformation $\cZ_G \Rightarrow \cZ^{\textrm{fin}}_G$, where $ \cZ^{\textrm{fin}}_G(M) =  \K{\Top/\Rep{G}(M)(\FF_q)}$ and analogously for bordisms. This TQFT computes the invariant $[\Rep{G}(M)(\FF_q)] \in \K{\Top}$, which is the number of $\FF_q$-points of $\Rep{G}(M)$. Moreover, an analogous argument to Proposition \ref{prop:TQFT-dim0-monoidal} shows that $\cZ^{\textrm{fin}}_{G}$ is a monoidal TQFT. Furthermore, the same construction for $\cZ^{\textrm{fin}}_G$ works verbatim in the case that $G$ is a finite group, giving rise to a monoidal TQFT that computes the number of points of the $G$-representation varieties.
\end{rmk}

We finish this section with a rather non-standard proof of a well-known result in group theory.

\begin{cor}
Let $G$ be a finite group with $c > 0$ conjugacy classes. Then $G$ has $|G|c$ pairwise commuting pairs of elements.

\begin{proof}
The twisting map $\cZ_{G}^{\textrm{fin}}(\Theta_{S^1}): (\K\Top)^{\oplus |G|} \to (\K\Top)^{\oplus |G|}$ is given by $\cZ_{G}^{\textrm{fin}}(\Theta_{S^1})(1_g)= \sum_{h} |\Stab(g)| 1_h$ where the sum runs over the elements $h \in G$ in the same conjugacy class as $g$. In particular, the trace of this map is 
$$
	\Tr(\cZ_{G}^{\textrm{fin}}(\Theta_{S^1})) = \sum_{\cC \subseteq G} |\Stab(\cC)| |\cC| = \sum_{\cC \subseteq G} |G| = |G|c,
$$
where the sum runs over the conjugacy classes $\cC$ of $G$. Since the number of pairwise commuting pairs of elements of $G$ is $\chi(S^1 \times S^1)/|G| = \Tr(\cZ_{G}^{\textrm{fin}}(\Theta_{S^1}^2))/|G| = \Tr(\cZ_{G}^{\textrm{fin}}(\Theta_{S_1}))$, the result follows.
\end{proof}
\end{cor}

\begin{rmk}
The previous count agrees with the usual counting formula for solutions of equations over finite fields. For instance, using \cite[Equation (2.3.8)]{Hausel-Rodriguez-Villegas:2008} we also get that the number of pairwise commuting elements of $G$ is
$$
	|\Hom(\ZZ \times \ZZ, G)| = |G| \sum_{\chi \in \textrm{Irr}(G)} \left(\frac{|G|}{\chi(1)}\right)^{2\cdot 1 - 2} = |G| \sum_{\chi \in \textrm{Irr}(G)}1 = |G|c,
$$
where the sum runs over the set of irreducible characters of $G$.
\end{rmk}

\subsection{$E$-polynomials of $\SL{2}(\CC)$-representation varieties as a dynamical system}\label{sec:sl2-rep-var}

In this section, we shall show how to apply the techniques from Section \ref{sec:dynamical-system} to provide a close formula for the $E$-polynomials of representation varieties for $G = \SL{2}(\CC)$. Even though the results of \cite{GP-2018} complete the proof in this case, for a general group $G$ the argument is somehow heuristic and does not lead to a rigorous proof. However, this idea can be used to exhibit conjectural closed formulas for the virtual classes of representation varieties for any group, which to our knowledge are currently unknown in many cases. 

First of all, notice that, in the case that $G$ is a complex group, an analogous construction to Theorem \ref{thm:lax-monoidal-rep-var} can be done to lax-quantize $E$-polynomials of representation varieties, leading to a functor
$$
	\cZ_G^{\textrm{Hod}}: \Bordp{n+1} \to \Mod{\K{\textbf{MHS}}},
$$
where $ \K{\textbf{MHS}}$ is the Grothendieck ring of mixed Hodge structures.
In this setting, instead of taking $\cZ_G(M,A) = \K{\Varrel{\Rep{G}(M,A)}}$, one should set $\cZ_G^{\textrm{Hod}}(M,A) = \K{\textbf{MHM}_{\Rep{G}(M,A)}}$, the Grothendieck ring of the category of Saito's mixed Hodge modules \cite{Saito:1990} over $\Rep{G}(M,A)$. The functor assigns to a closed manifold $(n+1)$-dimensional manifold $(W, A)$ the virtual class of the Hodge structure on the cohomology of the associated representation variety, $[H_c^\bullet(\Rep{G}(W,A))] \in \K{\textbf{MHS}}$. From this virtual Hodge structure, the $E$-polynomial is readily obtained by taking the weighted dimension of each Hodge piece. For the detailed construction of this TQFT, please refer to \cite{GPLM-2017,GP-2018}.

Observe that we can see $\ZZ[q]$ as a subring of $\K{\textbf{MHS}}$ by sending $q \mapsto [\QQ(-1)]$, the virtual class of the Tate motive of weight $2$. Since $\QQ(-1)$ is precisely the mixed Hodge structure of the affine line, this is compatible with both taking the $E$-polynomial and the natural map $\K{\Var{\CC}} \to \K{\textbf{MHS}}$. Further passing to the fraction field of $\ZZ[q]$, this gives rise to a lax monoidal TQFT
\begin{equation}\label{eq:TQFT-Epol}
	\hat{\cZ}_G^{\textrm{Hod}}: \Bordp{n+1} \to \Vect{\QQ(q)}.
\end{equation}
In particular, as pointed out in Remark \ref{rmk:almost-tqft-from-basepoints}, for surfaces this functor restricts to an almost-TQFT without basepoints $\hat{Z}_G^{\textrm{Hod}}: \Tub{1+1} \to \Vect{\QQ(q)}$. From the work of \cite{MM} (see also \cite[Section 5.4]{GP-2018}), an explicit expression for $\tilde{T} = \hat{Z}_{\SL{2}(\CC)}^{\textrm{Hod}}(L) :\hat{Z}_{\SL{2}(\CC)}^{\textrm{Hod}}(S^1, \star) \to \hat{Z}_{\SL{2}(\CC)}^{\textrm{Hod}}(S^1, \star)$ can be obtained, in particular showing that (\ref{eq:TQFT-Epol}) gives rise to an almost-TQFT taking values in finite dimensional vector spaces. The aim of this section is to show that, even without this knowledge, the map $\tilde{T}$ can be fully described from only a bunch of simple computations.

In fact, suppose that we have computed the virtual Hodge structures of the $\SL{2}(\CC)$-representation varieties over the genus $g$ surface $\Sigma_g$ for some finite number of genii $g$. This can be done through a long case-by-case analysis. The results for $0 \leq g \leq 11$ are shown in Table \ref{tab:E-pols}. 

\begin{table}[h!]
\centering
\tiny
 \begin{tabular}{|c | c|} 
 \hline
 \normalsize{$g$} & \normalsize{$[H_c^\bullet(\Rep{}(\Sigma_g))]$} \\ 
 \hline\hline
1 & $\begin{matrix*}[l]1\end{matrix*}$ \\\hline
2 & $\begin{matrix*}[l] q^{4} + 4  q^{3} - q^{2} - 4  q \end{matrix*}$ \\\hline
3 & $\begin{matrix*}[l] q^{9} + q^{8} + 12  q^{7} + 2  q^{6} - 3  q^{4} - 12  q^{3} - q \end{matrix*}$ \\\hline
4 & $\begin{matrix*}[l] q^{15} - 5  q^{13} + q^{12} + 73  q^{11} + 9  q^{10} + 295  q^{9} - 5  q^{8} - 295  q^{7} - 5  q^{6} - 73  q^{5} + 5  q^{3} - q \end{matrix*}$ \\\hline
5 & $\begin{matrix*}[l] q^{21} - 7  q^{19} + 21  q^{17} + q^{16} + 220  q^{15} + 20  q^{14} + 3584  q^{13} + 14  q^{12} - 28  q^{10} - 3584  q^{9} - 7  q^{8} - 220  q^{7} - 21  q^{5} + 7  q^{3} - q \end{matrix*}$ \\\hline
6 & $\begin{matrix*}[l] q^{27} - 9  q^{25} + 36  q^{23} - 84  q^{21} + q^{20} + 1149  q^{19} + 35  q^{18} + 27459  q^{17} + 90  q^{16}  + 43044  q^{15} - 42  q^{14} - 43044  q^{13} - 75  q^{12} \\\; - 27459  q^{11} - 9  q^{10} - 1149  q^{9} + 84  q^{7} - 36  q^{5} + 9  q^{3} - q \end{matrix*}$ \\\hline
7 & $\begin{matrix*}[l] q^{33} - 11  q^{31} + 55  q^{29} - 165  q^{27} + 330  q^{25} + q^{24} + 3633  q^{23} + 54  q^{22} + 180587  q^{21} + 275  q^{20} + 675235  q^{19} + 132  q^{18} - 297  q^{16} \\\; - 675235  q^{15} - 154  q^{14} - 180587  q^{13} - 11  q^{12} - 3633  q^{11} - 330  q^{9} + 165  q^{7} - 55  q^{5} + 11  q^{3} - q \end{matrix*}$ \\\hline
8 & $\begin{matrix*}[l] q^{39} - 13  q^{37} + 78  q^{35} - 286  q^{33} + 715  q^{31} - 1287  q^{29} + q^{28} + 18099  q^{27} + 77  q^{26} + 1063101  q^{25} + 637  q^{24} + 7029243  q^{23} \\ \; + 1001  q^{22} + 7026877  q^{21} - 429  q^{20} - 7026877  q^{19} - 1001  q^{18} - 7029243  q^{17} - 273  q^{16} - 1063101  q^{15} - 13  q^{14} - 18099  q^{13} \\\;+ 1287  q^{11} - 715  q^{9} + 286  q^{7} - 78  q^{5} + 13  q^{3} - q \end{matrix*}$ \\\hline
9 & $\begin{matrix*}[l] q^{45} - 15  q^{43} + 105  q^{41} - 455  q^{39} + 1365  q^{37} - 3003  q^{35} + 5005  q^{33} + q^{32} + 59100  q^{31} + 104  q^{30} + 5904480  q^{29} + 1260  q^{28} \\\;+ 59631040  q^{27} + 3640  q^{26} + 131201616  q^{25} + 1430  q^{24} - 3432  q^{22} - 131201616  q^{21} - 2548  q^{20} - 59631040  q^{19} - 440  q^{18} \\\;- 5904480  q^{17} - 15  q^{16} - 59100  q^{15} - 5005  q^{13} + 3003  q^{11} - 1365  q^{9} + 455  q^{7} - 105  q^{5} + 15  q^{3} - q \end{matrix*}$ \\\hline
10 & $\begin{matrix*}[l] q^{51} - 17  q^{49} + 136  q^{47} - 680  q^{45} + 2380  q^{43} - 6188  q^{41} + 12376  q^{39} - 19448  q^{37} + q^{36} + 286453  q^{35} + 135  q^{34} + 31170571  q^{33} \\\;+ 2244  q^{32} + 445660984  q^{31} + 9996  q^{30} + 1622121640  q^{29} + 11934  q^{28} + 1274542972  q^{27} - 4862  q^{26} - 1274542972  q^{25} \\\;- 13260  q^{24} - 1622121640  q^{23} - 5508  q^{22} - 445660984  q^{21} - 663  q^{20} - 31170571  q^{19} - 17  q^{18} - 286453  q^{17} + 19448  q^{15}  \\\;- 12376  q^{13}  + 6188  q^{11} - 2380  q^{9} + 680  q^{7} - 136  q^{5} + 17  q^{3} - q \end{matrix*}$ \\\hline
11 & $\begin{matrix*}[l] q^{57} - 19  q^{55} + 171  q^{53} - 969  q^{51} + 3876  q^{49} - 11628  q^{47} + 27132  q^{45} - 50388  q^{43} + 75582  q^{41} + q^{40} + 956197  q^{39}  + 170  q^{38} \\\; + 159475607  q^{37} + 3705  q^{36} + 3048129207  q^{35} + 23256  q^{34} + 16257140715  q^{33} + 48450  q^{32} + 26417750100  q^{31} + 16796  q^{30} \\\; - 41990  q^{28} - 26417750100  q^{27} - 38760  q^{26} - 16257140715  q^{25} - 10659  q^{24} - 3048129207  q^{23} - 950  q^{22} - 159475607  q^{21} \\\; - 19  q^{20} - 956197  q^{19} - 75582  q^{17} + 50388  q^{15} - 27132  q^{13} + 11628  q^{11} - 3876  q^{9} + 969  q^{7} - 171  q^{5} + 19  q^{3} - q \end{matrix*}$ \\\hline
 \end{tabular}
 \vspace{0.2cm}
 \caption{Virtual Hodge structures of $\SL{2}(\CC)$-representation varieties for closed orientable surfaces of genus $1 \leq g \leq 11$. Here, we set $q = [H^\bullet_c(\CC)]$ the virtual Hodge structure of the affine line.}\label{tab:E-pols}
\end{table}

Following Remark \ref{rmk:computation-tqft}, let us form the vectors
\begin{align}\label{eq:vector-rec}
	v_{g, N} = \left([H_c^\bullet(\Rep{}(\Sigma_g))], \ldots, [H_c^\bullet(\Rep{}(\Sigma_{g+N-1}))] \right) \in \ZZ[q]^N \subseteq \K{\textbf{MHS}}^N.
\end{align}

A direct calculation shows that for $N \leq 6$, the vectors $v_{0, N}, \ldots, v_{N-1,N}$ are linearly independent over $\QQ(q)$ but for $N=7$ they are linearly dependent. This suggests that the dimension of the underlying vector space of the almost-Frobenius algebra is $N = 6$. A direct computation shows that
$$
	v_{6,6} = P_0 v_{0,6} + P_1 v_{1,6} + \ldots + P_5 v_{5,6},
$$
where the polynomials $P_i \in \ZZ[q]$ are
\begin{align*}
P_0 &= q^{6} + 9  q^{4} + 9  q^{2} + 1, \\
P_1 &= -11  q^{10} - 29  q^{8} + 16  q^{6} - 29  q^{4} - 11  q^{2}, \\
P_2 &= 43  q^{14} - 25  q^{12} - 18  q^{10} - 18  q^{8} - 25  q^{6} + 43  q^{4}, \\
P_3 &= -73  q^{18} + 198  q^{16} - 135  q^{14} + 20  q^{12} - 135  q^{10} + 198  q^{8} - 73  q^{6}, \\
P_4 &= 56  q^{22} - 280  q^{20} + 504  q^{18} - 280  q^{16} - 280  q^{14} + 504  q^{12} - 280  q^{10} + 56  q^{8}, \\
P_5 &= -16  q^{26} + 128  q^{24} - 448  q^{22} + 896  q^{20} - 1120  q^{18} + 896  q^{16} - 448  q^{14} + 128  q^{12} - 16  q^{10}. \\
\end{align*}
As argued in Remark \ref{rmk:computation-tqft}, this implies the recurrence relation
$$
	\left[H_c^\bullet(\Rep{}(\Sigma_{g}))\right] = P_0 \left[H_c^\bullet(\Rep{}(\Sigma_{g-6}))\right] + P_1 \left[H_c^\bullet(\Rep{}(\Sigma_{g-5}))\right] + \ldots + P_5 \left[H_c^\bullet(\Rep{}(\Sigma_{g-1}))\right],
$$
for all $g \geq 6$.

In other words, we have that the map $\tilde{T}: \QQ(q)^6 \to \QQ(q)^6$ of the almost-Frobenius algebra is the only map satisfying $\tilde{T}(v_{g, 6}) = v_{g+1, 6}$ for $g < 5$ and $\tilde{T}(v_{5, 6}) = P_0 v_{0,6} + \ldots + P_5 v_{5,6}$. That is, in the basis $v_{0, 6}, \ldots, v_{5,6}$, we have

$$
	\tilde{T} = \left(\begin{matrix} 0 & 0 & 0 & 0 & 0 & P_0 \\
	1 & 0 & 0 & 0 & 0 & P_1 \\
	0 & 1 & 0 & 0 & 0 & P_2 \\
	0 & 0 & 1 & 0 & 0 & P_3 \\
	0 & 0 & 0 & 1 & 0 & P_4 \\
	0 & 0 & 0 & 0 & 1 & P_5 \\
	 \end{matrix}\right).
$$

Moreover, setting $Q_0 = \left(\begin{array}{c|c|c|c|c|c} v_{0,6} & v_{1,6} & v_{2,6} & v_{3,6} & v_{4,6} & v_{5,6} \end{array}\right)$, we have that $Q_0\tilde{T}Q_0^{-1}$ is the matrix of $\tilde{T}$ in the canonical basis. In this basis, we recover the whole almost-Frobenius algebra associated to the lax monoidal TQFT of Theorem \ref{thm:lax-monoidal-rep-var} with the inclusion $\tilde{\epsilon}(1) = v_{0,n}$ and the projection onto the first component $\tilde{\eta}: \QQ(q)^6 \to \QQ(q)$. For instance, writing $\tilde{T} = QDQ^{-1}$ with $D$ a diagonal matrix and setting $e_1 = (1,0,0,0,0,0)$, we get that
\begin{align}\label{for:motive-SL2}
\begin{split}
	[H_c^\bullet(\Rep{}(\Sigma_g))] &= \tilde{\eta}(\tilde{T}^g(v_{0,6})) = e_1^t Q_0QD^gQ^{-1}Q_0^{-1} e_1 \\
	&= \,{\left(q^2 - 1\right)}^{2g - 1} q^{2g - 1} +
\frac{1}{2} \, {\left(q -
1\right)}^{2g - 1}q^{2g -
1}(q+1){\left({2^{2g} + q - 3}\right)} 
\\ & \qquad+ \frac{1}{2} \,
{\left(q + 1\right)}^{2g + r - 1} q^{2g - 1} (q-1){\left({2^{2g} +q -1}\right)} + q(q^2-1)^{2g-1}.
\end{split}
\end{align}
This formula agrees with \cite[Proposition 11]{MM} (see also \cite[Remark 5.11]{GP-2018}).

\begin{rmk}
The map $\tilde{T}$ computed above does not coincide exactly with $\hat{Z}_{\SL{2}(\CC)}^{\textrm{Hod}}(L)$. The reason is that $\tilde{T}^g$ is computing $[H_c^\bullet(\Rep{\SL{2}(\CC)}(\Sigma_g))]$ but, with the notation of Remark \ref{rmk:almost-tqft-from-basepoints}, $\hat{Z}_{\SL{2}(\CC)}^{\textrm{Hod}}(L^g) = \hat{\cZ}_{\SL{2}(\CC)}^{\textrm{Hod}}(\kappa(L^g))$ is computing the virtual Hodge structure $[H_c^\bullet(\SL{2}(\CC))]^{g}[\Rep{\SL{2}(\CC)}(\Sigma_g)]$, since $\kappa(L^g)$ is equipped with $g+1$ basepoints. In other words, in this section we have actually computed $\tilde{T} = \frac{1}{[H_c^\bullet(\SL{2}(\CC))]}\hat{Z}_{\SL{2}(\CC)}^{\textrm{Hod}}(L)$.
\end{rmk}

\begin{rmk}
The previous argument is not completely rigorous since, with this method, we cannot prove that the predicted almost-TQFT of dimension $6$ actually exists. We found the predicted dimension by a na\"ive linear independence argument, but we do not know whether the lax monoidal TQFT of (\ref{eq:TQFT-Epol}) can be restricted to an almost-TQFT with values in finite dimensional vector spaces, since we only tested it on finitely many cases. This is the only gap in the argument. However, it turns out that such property was proven in \cite[Theorem 4.9]{GP-2018} for $G = \SL{2}(\CC)$, so formula (\ref{for:motive-SL2}) actually holds true.
\end{rmk}

\begin{rmk}
Indeed, this property of having a finite dimensional invariant subspace cannot be detected through this kind of arguments. For instance, let $e_i$ be the $i$-th vector of the canonical basis and set
$$
	v_1 = e_1, \quad v_2 = e_2, \quad v_g = v_{g-2} + v_{g-1} + e_{g} \textrm{  for }g \geq 3.
$$
Then, for $V = \langle v_g \rangle_{g = 1}^\infty$ and the map $\tilde{T}: V \to V$ given by $\tilde{T}(v_g) = v_{g+1}$ there is no such invariant subspace. However, if we take $\tilde{\eta}(v_1) = \tilde{\eta}(v_2) = 1$ and $\tilde{\eta}(e_g) = 0$ for $g \geq 3$, we have that $\tilde{\eta}(v_g)$ is the $g$-th Fibonacci number for $g \geq 1$. Hence, the vectors (\ref{eq:vector-rec}) form a linear recurrence of order two but they cannot detect that there exists no finitely generated invariant submodule.
\end{rmk}

The upshot of this computation is that from the simple knowledge of the virtual Hodge class of the representation variety for only finitely many genii, which can be computed through a brute force approach (see for instance the algorithm described in \cite[Appendix A]{hablicsek2022virtual} or the algebraic representatives method developed in \cite{vogel2023motivic}), the virtual Hodge class for arbitrary genus can be computed thanks to the existence of an underlying TQFT.

%\nocite{*}

\bibliography{bibliography.bib}{}
\bibliographystyle{abbrv}

\end{document}